\documentclass[11pt]{article}
\usepackage[margin=1in]{geometry} 
\usepackage{mathtools}

\usepackage[LGR,T1]{fontenc}
\DeclareSymbolFont{upgreek}{LGR}{cmr}{m}{n}
\DeclareMathSymbol{\deltaup}{\mathord}{upgreek}{`d}
\DeclareMathSymbol{\piup}{\mathord}{upgreek}{`p}
\DeclareMathSymbol{\epsilonup}{\mathord}{upgreek}{`e}

\DeclareMathSymbol{\omegaup}{\mathord}{upgreek}{`w}

\newcommand{\dynkinradius}{.075cm}
\newcommand{\dynkinstep}{1cm}
\newcommand{\dynkindot}[2]{\fill (\dynkinstep*#1,\dynkinstep*#2) circle (\dynkinradius);}

\newcommand{\dynkinline}[4]{\draw[thin] (\dynkinstep*#1,\dynkinstep*#2) -- (\dynkinstep*#3,\dynkinstep*#4);}
\newcommand{\dynkindots}[4]{\draw[dotted] (\dynkinstep*#1,\dynkinstep*#2) -- (\dynkinstep*#3,\dynkinstep*#4);}
\newcommand{\dynkindoubleline}[4]{\draw[double,postaction={decorate}] (\dynkinstep*#1,\dynkinstep*#2) -- (\dynkinstep*#3,\dynkinstep*#4);}

\newenvironment{dynkin}{\begin{tikzpicture}[decoration={markings,mark=at position 0.7 with {\arrow{>}}}]}
{\end{tikzpicture}}

\makeatletter
\newcommand*{\shifttext}[2]{%
  \settowidth{\@tempdima}{#2}%
  \makebox[\@tempdima]{\hspace*{#1}#2}%
}
\makeatother

\newcommand{\uniondot}{{\mathrlap{\raisebox{0.4pt}{\shifttext{5pt}{$\cdot$}}}{\cup}}}

\usepackage{amssymb,amsfonts,amsmath,bbm,mathrsfs,stmaryrd}
\usepackage{soul,xcolor}
\usepackage{url}
\usepackage{ytableau}

\usepackage[colorlinks,
             linkcolor=black!75!red,
             citecolor=blue,
             pdftitle={Root-Reductive},
             pdfauthor={Thanasin Nampaisarn},
             pdfproducer={pdfLaTeX},
             pdfpagemode=None,
             bookmarksopen=true
             bookmarksnumbered=true]{hyperref}

\usepackage{xr}

\usepackage{tikz,tikz-cd}
\usetikzlibrary{matrix,arrows,calc,decorations.pathreplacing,decorations.markings,intersections,shapes.geometric,through,fit,shapes.symbols,positioning,decorations.pathmorphing,snakes}
\everymath{\displaystyle}

\usepackage[amsmath,thmmarks,hyperref]{ntheorem}
\usepackage{cleveref}

\crefname{section}{Section}{Sections}
\crefformat{section}{#2Section~#1#3} 
\Crefformat{section}{#2Section~#1#3} 

\crefname{subsection}{\S}{\S\S}
\crefformat{subsection}{#2\S#1#3} 
\Crefformat{subsection}{#2\S#1#3} 

\theoremstyle{plain}

\newtheorem{lem}{Lemma}[section]
\newtheorem{prop}[lem]{Proposition}
\newtheorem{cor}[lem]{Corollary}
\newtheorem{thm}[lem]{Theorem}

\newtheorem{openq}[lem]{Open Question}

\theoremstyle{nonumberplain}

\newtheorem{qf'bis}{\Cref{le.quot_filt'} bis}

\theoremstyle{plain}
\theorembodyfont{\upshape}
\theoremsymbol{\ensuremath{\blacklozenge}}

\newtheorem{define}[lem]{Definition}

\crefname{definition}{definition}{definitions}
\crefformat{definition}{#2definition~#1#3} 
\Crefformat{definition}{#2Definition~#1#3} 

\crefname{ex}{example}{examples}
\crefformat{example}{#2example~#1#3} 
\Crefformat{example}{#2Example~#1#3} 

\crefname{remark}{remark}{remarks}
\crefformat{remark}{#2remark~#1#3} 
\Crefformat{remark}{#2Remark~#1#3} 

\crefname{convention}{convention}{conventions}
\crefformat{convention}{#2convention~#1#3} 
\Crefformat{convention}{#2Convention~#1#3}

\crefname{lemma}{lemma}{lemmas}
\crefformat{lemma}{#2lemma~#1#3} 
\Crefformat{lemma}{#2Lemma~#1#3} 

\crefname{proposition}{proposition}{propositions}
\crefformat{proposition}{#2proposition~#1#3} 
\Crefformat{proposition}{#2Proposition~#1#3} 

\crefname{corollary}{corollary}{corollaries}
\crefformat{corollary}{#2corollary~#1#3} 
\Crefformat{corollary}{#2Corollary~#1#3} 

\crefname{theorem}{theorem}{theorems}
\crefformat{theorem}{#2theorem~#1#3} 
\Crefformat{theorem}{#2Theorem~#1#3} 

\crefname{assumption}{assumption}{Assumptions}
\crefformat{assumption}{#2assumption~#1#3} 
\Crefformat{assumption}{#2Assumption~#1#3} 

\crefname{equation}{}{}
\crefformat{equation}{(#2#1#3)} 
\Crefformat{equation}{(#2#1#3)}

\theoremstyle{nonumberplain}
\theoremsymbol{\ensuremath{\blacksquare}}

\newtheorem{pf}{Proof}

\newtheorem{proof-of-uptri}{Proof of \Cref{pr.up-tri}}
\newtheorem{proof-of-simples}{Proof of \Cref{pr.simples}}
\newtheorem{proof-of-locnoe}{Proof of \Cref{pr.loc-noe}}
\newtheorem{proof-of-exts}{Proof of \Cref{th.exts}}
\newtheorem{proof-of-univ}{Proof of \Cref{th.univ}}
\newtheorem{proof-of-univ-gen}{Proof of \Cref{th.univ-gen}}

\newcommand{\tendsto}{\boldsymbol{\to}}

\newcommand{\qedhere}{\mbox{}\hfill\ensuremath{\blacksquare}}

\usepackage{graphics}
\usepackage{epsfig}
\usepackage{epstopdf}
\usepackage{slashed}
\usepackage{enumerate}
\usepackage{scalerel}
\usepackage{epic}
\usepackage{indentfirst}
\usepackage{verbatim}

\usepackage{etoolbox}

\usepackage[T1]{fontenc}
\usepackage{calligra}
\usepackage{xcolor}
\makeatletter
\let\ps@plainorig\ps@plain
\makeatother

\newcommand\swabfamily{\usefont{U}{yswab}{m}{n}}
\DeclareTextFontCommand{\textswab}{\swabfamily}
\newcommand\frakfamily{\usefont{U}{yfrak}{m}{n}}
\DeclareTextFontCommand{\textfrak}{\frakfamily}
\newcommand\gothfamily{\usefont{U}{ygoth}{m}{n}}
\DeclareTextFontCommand{\textgoth}{\gothfamily}

\DeclareMathOperator{\Ext}{Ext}

\DeclareMathOperator{\Hom}{Hom}
\DeclareMathOperator{\End}{End}
\DeclareMathOperator{\im}{im}
\DeclareMathOperator{\id}{id}
\DeclareMathOperator{\ch}{ch}
\DeclareMathOperator{\Res}{Res}

\DeclareFontFamily{U}{mathb}{\hyphenchar\font45}
\DeclareFontShape{U}{mathb}{m}{n}{
      <5> <6> <7> <8> <9> <10> gen * mathb
      <10.95> mathb 10 <12> <14.4> <17.28> <20.74> <24.88> mathb12
      }{}
\DeclareSymbolFont{mathb}{U}{mathb}{m}{n}
\DeclareMathSymbol{\precneq}{3}{mathb}{"AC}

\DeclareSymbolFont{euletters}{U}{eur}{m}{n}
\SetSymbolFont{euletters}{bold}{U}{eur}{b}{n}
\DeclareMathSymbol{\varp}{\mathalpha}{euletters}{"7D}

\newcommand{\bggO}{\pazocal{O}}
\newcommand{\bbggO}{{\bar{\pazocal{O}}}}

\newcommand{\defeq}{\overset{\mathrm{def}}{=\!=}}

\newcommand{\amsbb}[1]{\mathbb{#1}}

\def\suchthat#1{%
\,\pmb{\left|\vphantom{#1}\right.}\,%
#1
}

\newcommand{\gfrak}{\mathfrak{g}}
\newcommand{\hfrak}{\mathfrak{h}}
\newcommand{\kfrak}{\mathfrak{k}}
\newcommand{\bfrak}{\mathfrak{b}}
\newcommand{\nfrak}{\mathfrak{n}}
\newcommand{\pfrak}{\mathfrak{p}}

\newcommand{\gllie}{\textup{\swabfamily gl}}
\newcommand{\sllie}{\textup{\swabfamily s:l}}
\newcommand{\solie}{\textup{\swabfamily s:o}}
\newcommand{\splie}{\textup{\swabfamily s:p}}

\newcommand{\Ulie}{\textup{\swabfamily U}}
\newcommand{\Mlie}{\textup{\swabfamily M}}
\newcommand{\Nlie}{\textup{\swabfamily N}}
\newcommand{\Llie}{\textup{\swabfamily L}}
\newcommand{\Plie}{\textup{\swabfamily P}}

\newcommand{\Ilie}{\textup{\swabfamily I}}

\newcommand{\Vlie}{\textup{\swabfamily V}}

\newcommand{\tcyr}{\textup{\textsc{t}}}
\usepackage{calrsfs}
\DeclareMathAlphabet{\pazocal}{OMS}{zplm}{m}{n}

\title{On Categories $\bggO$ for Root-Reductive Lie Algebras}
\author{Thanasin Nampaisarn}

\begin{document}
\setstcolor{red}

\date{}

\newcommand{\Addresses}{{
  \bigskip
  \footnotesize

  \textsc{Jacobs University Bremen, Campus Ring 1, 28759 Bremen, Germany}\par\nopagebreak
  \textit{E-mail address}: \texttt{t.nampaisarn@jacobs-university.de}
}}

\maketitle

\begin{abstract}
Let $\gfrak$ be a root-reductive Lie algebra over an algebraically closed field $\amsbb{K}$ of characteristic $0$ with a splitting Borel subalgebra $\bfrak$ containing a splitting maximal toral subalgebra $\hfrak$.  We study the category $\bbggO$ consisting of all $\hfrak$-weight $\gfrak$-modules which are locally $\bfrak$-finite and have finite-dimensional $\hfrak$-weight spaces.  The focus is on very special Borel subalgebras called the Dynkin Borel subalgebras.  This paper serves as an initial passage to the understanding of categories $\bggO$ for infinite-dimensional root-reductive Lie algebras.
\end{abstract}

\noindent {\em Key words: root-reductive Lie algebras, finitary Lie algebras, highest-weight modules, BGG categories $\bggO$, truncated subcategories}

\vspace{.5cm}

\noindent{MSC 2010: 17B10; 17B20; 17B22}


\section*{Introduction}
\addcontentsline{toc}{section}{Introduction}

The purpose of this paper is to introduce an analogue of the Bernstein-Gel'fand-Gel'fand category $\bggO$ for a class of infinite-dimensional Lie algebras over an algebraically closed field $\amsbb{K}$ of characteristic $0$: the root-reductive Lie algebras.  These Lie algebras arise as direct limits of reductive Lie algebras with certain restrictions.  Examples of such Lie algebras are the classical finitary Lie algebras $\gllie_\infty(\amsbb{K})$, $\sllie_\infty(\amsbb{K})$, $\solie_\infty(\amsbb{K})$, and $\splie_\infty(\amsbb{K})$, which are well understood (see \cite{BS} and \cite{penkov-strade}).

	 Extensive studies of analogues of categories $\bggO$ have also been undertaken for affine Kac-Moody Lie algebras (see \cite{FF} and \cite{SVV}).  In a parallel fashion to the finite-dimensional theory, an affine Kac-Moody Lie algebra $\hat{\gfrak}$ has a Borel subalgebra $\hat{\bfrak}$ containing a Cartan subalgebra $\hat{\hfrak}$.  It is natural to define the category $\bggO$ for the Lie algebra $\hat\gfrak$ by requiring that $\hat{\hfrak}$ act semisimply and $\hat{\bfrak}$ act locally finitely on each module in $\bggO$.  
	 
	 Similarly to the affine Kac-Moody case, splitting maximal toral subalgebras and splitting Borel subalgebras play essential roles in our approach.  Both maximal toral subalgebras and Borel subalgebras of root-reductive Lie algebras have been studied in \cite{dancohen}, \cite{BorelDC}, \cite{DPS}, and \cite{penkov-dimitrov}.  Furthermore, \cite{DPS} proves that, for every root-reductive Lie algebra $\gfrak$, the derived ideal $[\gfrak,\gfrak]$ is a direct sum of the finite-dimensional simple Lie algebras and the simple finitary Lie algebras $\sllie_\infty$, $\solie_\infty$, $\splie_\infty$, each occuring with at most countable multiplicity.
	 
	  To define a category $\bggO$ for a root-reductive Lie algebra $\gfrak$, we further need to understand the structure of Borel subalgebras of $\gfrak$.  Borel subalgebras for the simple finitary Lie algebras $\sllie_\infty$, $\solie_\infty$, and $\splie_\infty$ and for root-reductive Lie algebras are very well understood (see~\cite{BorelDC} and~\cite{DP1999}).  For a given Borel subalgebra $\bfrak$ of $\gfrak$ containing a splitting maximal toral subalgebra $\hfrak$, we define an extended category $\bggO$, denoted by $\bbggO$, for $\gfrak$ with respect to $\bfrak$ by demanding that the objects in the category be $\hfrak$-semisimple with locally finite $\bfrak$-action and with finite-dimensional $\hfrak$-weight spaces.

	While the category $\bbggO$ has its own merits for an arbitrary splitting Borel subalgebra $\bfrak$, a special class of Borel subalgebras, known as the Dynkin Borel subalgebras, eases the study of the category $\bbggO$.  The category $\bbggO$ with respect to a Dynkin Borel subalgebra $\bfrak$ has many additional desirable properties.  For example, it contains all Verma modules, every object in $\bbggO$ has an analogue of composition series, and the Kazhdan-Lusztig theory generalizes to the category $\bbggO$.

This paper consists of four sections.  The first section introduces the root-reductive Lie algebras, important subalgebras such as splitting maximal toral subalgebras and
splitting Borel subalgebras, the simple finitary Lie algebras, and Verma modules.  The second section defines the extended category $\bggO$ with respect to a Dynkin Borel subalgebra of a root-reductive Lie algebra and studies the general properties of objects in this category.  The third section explores the block structure of $\bbggO$ as well as the related Kazhdan-Lusztig theory.  The final section deals with the obstruction within $\bbggO$---the lack of injective objects.  This fourth section follows the idea of~\cite{RCW} to establish a version of BGG reciprocity of the truncated subcategories of $\bbggO$.

\subsection*{Acknowledgement}

The author has been supported in part by DFG (Deutsche Forschungsgemeinschaft) via the priority program
Representation Theory.  The author would also like to thank Professor Ivan Penkov from Jacobs University Bremen for his tremendous help towards the completion of this paper.

\section{Preliminaries}
\label{ch:intro}

The base field is $\amsbb{K}$, which is assumed to be algebraically closed and of characteristic $0$.  All Lie algebras and vector spaces are defined over $\amsbb{K}$.   For a vector space $V$, $V^*$ denotes its algebraic dual space $\Hom_\amsbb{K}(V,\amsbb{K})$.

\subsection{Root-Reductive Lie Algebras}
\label{ch:localfiniteness}

	A Lie algebra $\mathfrak{g}$ is said to be \emph{locally finite}\index{local finiteness} if every finite subset of $\mathfrak{g}$ generates a finite-dimensional subalgebra of $\mathfrak{g}$.\index{subalgebra}\index{subalgebra of $\gfrak$|see {subalgebra}}  A locally finite Lie algebra $\mathfrak{g}$ is \emph{locally solvable}\index{local solvability} if every finite-dimensional subalgebra of $\mathfrak{g}$ is solvable.  Similarly, a locally finite Lie algebra $\mathfrak{g}$ is \emph{locally nilpotent}\index{local nilpotency} if every finite-dimensional subalgebra of $\mathfrak{g}$ is nilpotent.\index{subalgebra of $\gfrak$|see {subalgebra}}  


\begin{define}
	An inclusion of finite-dimensional reductive Lie algebras $\mathfrak{l}\hookrightarrow\mathfrak{m}$ is a \emph{root inclusion} if, for some Cartan subalgebra $\mathfrak{c}$ of $\mathfrak{m}$, the subalgebra $\mathfrak{l}\cap \mathfrak{c}$ is a Cartan subalgebra of $\mathfrak{l}$ and each $(\mathfrak{l}\cap\mathfrak{c})$-root space  of $\mathfrak{l}$ is a root space of $\mathfrak{m}$.
	
	A \emph{root-reductive} Lie algebra $\gfrak$ is a locally finite Lie algebra $\gfrak=\bigcup_{n\in\amsbb{Z}_{> 0}}\,\gfrak_n$, where $\left(\gfrak_n\right)_{n\in\amsbb{Z}_{>0}}$ is a nested system of finite-dimensional reductive subalgebras, with the property that there exist nested subalgebras $\kfrak_1\subseteq \kfrak_2\subseteq \ldots$, where $\kfrak_n\subseteq \gfrak_n$ is a Cartan subalgebra of $\gfrak_n$, such that each inclusion $\gfrak_n \hookrightarrow \gfrak_{n+1}$ is a root inclusion with respect to the Cartan subalgebra $\kfrak_{n+1}$ of $\gfrak_{n+1}$.
\end{define}

\begin{define}
	Let $\gfrak$ be a root-reductive Lie algebra.  A subalgebra $\hfrak$ of $\gfrak$ is said to be a \emph{splitting maximal toral subalgebra}\index{subalgebra!splitting maximal toral subalgebra} if there exists a directed system $\left(\gfrak_n\right)_{n\in \amsbb{Z}_{>0}}$ of finite-dimensional reductive subalgebras of $\gfrak$ for which $\lim\limits_{\underset{n}{\boldsymbol{\longrightarrow}}}\, \gfrak_n=\gfrak$, $\hfrak\cap\gfrak_n$ is a maximal toral subalgebra of $\gfrak_n$ for each $n \in \amsbb{Z}_{>0}$, and $\gfrak$ has the \emph{$\hfrak$-root space decomposition} $
		\mathfrak{g}=\bigoplus\limits_{\alpha\in\mathfrak{h}^*}\mathfrak{g}^\alpha=\mathfrak{h}\,\oplus\,\bigoplus\limits_{\alpha \in \Delta}\mathfrak{g}^\alpha$,
		where $\mathfrak{g}^\alpha$ is the eigenspace\index{root space} $\big\{x\in\mathfrak{g}\suchthat{[h,x]=\alpha(h)x}\textup{ for all }h\in\hfrak\big\}$ for every $\alpha\in\mathfrak{h}^*$, and $\Delta$ is the set of \emph{roots}\index{root}, i.e.,  the nonzero linear functionals $\alpha \in \mathfrak{h}^*$ for which $\mathfrak{g}^\alpha \neq 0$.  For $\alpha \in \Delta$, $\gfrak^\alpha$ is known as the \emph{root space} associated to $\alpha$.   Note that $\gfrak^\alpha$ is always a one-dimensional vector space.
\end{define}

\begin{define}
	Let $\gfrak$ be a root-reductive Lie algebra.  A subalgebra $\bfrak$ of $\gfrak$ is said to be a \emph{splitting Borel subalgebra}\index{subalgebra!splitting Borel subalgebra} if $\bfrak$ is a maximal locally solvable subalgebra of $\gfrak$ containing a splitting maximal toral subalgebra of $\gfrak$.
\end{define}
	
	In~\cite{dancohen},~\cite{BorelDC},~\cite{DPS} and~\cite{DP1999}, root-reductive Lie algebras are studied.  In the case of root-reductive Lie algebra, a (splitting) Borel subalgebra $\bfrak$ containing a splitting maximal toral subalgebra $\hfrak$ arises from a choice of a set of \emph{$\bfrak$-positive roots} $\Delta^+\subseteq \Delta$ subject to the requirements that both subsets $\Delta^+$ and $\Delta\smallsetminus\Delta^+$ are additive and that $\alpha\in\Delta^+$ if and only if $-\alpha \in \Delta\smallsetminus\Delta^+$.  We set $\Delta^-:=-\Delta^+=\Delta\smallsetminus\Delta^+$ and call $\Delta^-$ the set of \emph{$\bfrak$-negative roots}.  Then $\Delta$ is the disjoint union $\Delta^+\uniondot \Delta^-$, the locally nilpotent subalgebra $\nfrak=\nfrak^+:=[\bfrak,\bfrak]$ is the direct sum $\bigoplus_{\alpha \in \Delta^+}\,\gfrak^\alpha$, and the Borel subalgebra $\bfrak$ is given by $\bfrak=\bfrak^+=\hfrak\oplus\nfrak^+$ (this is a direct sum of vector spaces, not of Lie algebras).  The Lie algebra $\gfrak$ has the \emph{triangular decomposition}
	$
		\gfrak=\nfrak^-\oplus\hfrak\oplus\nfrak^+$,
	where $\nfrak^-$ is the opposite subalgebra to $\nfrak^+$, namely, $\nfrak^-=\bigoplus_{\alpha\in\Delta^-}\,\gfrak^\alpha$.  The Lie algebra $\hfrak\oplus\nfrak^-$ is denoted by $\bfrak^-$.  It is opposite to $\bfrak^+$ in the sense that $\bfrak^++\bfrak^-=\gfrak$ and $\bfrak^+\cap\bfrak^-=\hfrak$.  By the  Poincar\'e-Birkhoff-Witt (PBW) Theorem, we see that\footnote{Here, $\Ulie(L)$ is the universal enveloping algebra of a Lie algebra $L$.} $\Ulie(\mathfrak{g})=\Ulie\left(\mathfrak{n}^-\right)\boldsymbol{\cdot} \Ulie\left(\mathfrak{h}\right)\boldsymbol{\cdot}\Ulie\left(\mathfrak{n}^+\right)$.  
	
	Let the \emph{root lattice} $\Lambda$ be the $\amsbb{Z}$-span in $\hfrak^*$ of $\Delta$, and $\Lambda^+$ be the $\amsbb{Z}_{\geq 0}$-span in $\hfrak^*$ of $\Delta^+$.  We equip $\mathfrak{h}^*$ with a partial order $\preceq$ given as follows:
		$
			\lambda \preceq \mu \textup{ iff } \mu- \lambda \in \Lambda^+
		$ 		
	 for all $\lambda,\mu\in\mathfrak{h}^*$.  We also write $\Lambda^-:=-\Lambda^+$ for the $\amsbb{Z}_{\geq 0}$-span of $\Delta^-$.

	\begin{define}
		An element $\alpha \in \Delta^+$ is said to be a \emph{simple $\bfrak$-positive root}, or a \emph{simple root} with respect to $\bfrak$, if it cannot be decomposed as a (finite) sum of two or more $\bfrak$-positive roots.  We usually use the symbol $\Sigma^+$ or $\Sigma$ for the set of simple $\bfrak$-positive roots. Similarly, we say that $\alpha \in \Delta^-$ is a \emph{simple negative root} with respect to $\bfrak$ if $-\alpha$ is a simple positive root.
	\end{define}
	
	From now on, $\gfrak$ is a root-reductive Lie algebra with a splitting Borel subalgebra $\bfrak$ containing a splitting maximal toral subalgebra $\hfrak$.  We assume that $\gfrak$ is the union $\bigcup_{n\in\amsbb{Z}_{>0}}\,\gfrak_n$ of nested finite-dimensional reductive subalgebras $\gfrak_n$ for which $\hfrak_n:=\hfrak\cap \gfrak_n$ is a maximal toral subalgebra of $\gfrak_n$, $\bfrak_n=\bfrak^+_n:=\bfrak\cap \gfrak_n$ is a Borel subalgebra of $\gfrak_n$, and $\nfrak_n=\nfrak^+_n:=\nfrak^+\cap\gfrak_n=\left[\bfrak^+_n,\bfrak^+_n\right]$ is a nilpotent subalgebra of $\gfrak_n$.   We also write $\bfrak_n^-:=\bfrak^-\cap\gfrak_n$ and $\nfrak^-_n:=\nfrak^-\cap\gfrak_n$.  In the case where $\gfrak$ is locally semisimple, we also assume that each $\gfrak_n$ is semisimple.  
	
	For each $n\in\amsbb{Z}_{>0}$, $W_n$ denotes the Weyl group of $\gfrak_n$.  Since the embedding $\gfrak_n\hookrightarrow \gfrak_{n+1}$ is a root inclusion, this induces an embedding $W_n\hookrightarrow W_{n+1}$.  The Weyl group $W$ of $\gfrak$ is simply the direct limit $\lim\limits_{\underset{n}{\boldsymbol{\longrightarrow}}}\,W_n$.

\subsection{Dynkin Borel Subalgebras}

In this subsection, $\gfrak$ need not be finite-dimensional.  Furthermore, let $\rho_n$ denote the half sum of positive roots of $\gfrak_n$ with respect to $\bfrak_n$ (sometimes, the linear functionals $\rho_n$ are called the \emph{local half sum of positive roots}).
	
	\begin{define}
		We say that $\bfrak$ is a \emph{Dynkin Borel subalgebra}  of $\gfrak$ if it is generated by $\hfrak$ and the simple root spaces.  
	\end{define}
	
	\begin{define}
		A $\bfrak$-positive root $\alpha$ is \emph{of finite length}  (with respect to $\bfrak$) if there are only finitely many ways to express it as a sum of positive roots; otherwise, $\alpha$ is \emph{of infinite length}  (with respect to $\bfrak$).  A $\bfrak$-negative root $\alpha$ is said to be \emph{of finite length}  (with respect to $\bfrak$) if the positive root $-\alpha$ is of finite length; otherwise, $\alpha$ is \emph{of infinite length}  (with respect to $\bfrak$). 
	\end{define}
	
	It is an easy exercise to prove that $\bfrak$ is Dynkin if and only if every root is of finite length (with respect to $\bfrak$).   In other words, a Dynkin Borel subalgebra $\bfrak$ is the direct sum of $\hfrak$ and the root spaces corresponding to roots of finite length.  
	
	The following proposition and its corollary are essential in this paper.  The proofs are straightforward, and therefore omitted.
	
	\begin{prop}
		Let $\bfrak$ be a Dynkin borel subalgebra of $\gfrak$.  Then,
		$
			\left.\rho_{n+1}\right|_{\hfrak_n}=\rho_n
		$
		for every positive integer $n$.
		\label{prop:hspr}
	\end{prop}
	
	\begin{cor}
		There exists $\rho\in\hfrak^*$ such that $\rho|_{\hfrak_n}=\rho_n$ for every $n\in\amsbb{Z}_{> 0}$ if and only if $\bfrak$ is a Dynkin Borel subalgebra.  That is, a \emph{global half sum of $\bfrak$-positive roots} $\rho$ is well defined.  Furthermore, if $\gfrak$ is locally semisimple (i.e., each $\gfrak_n$ is semisimple) and $\bfrak$ is Dynkin, then $\rho$ is unique.
		\label{cor:hspr}
	\end{cor}

	From the corollary above, Dynkin Borel subalgebras play a distinguished role because the existence of the global half sum of positive roots $\rho$ allows us to define the \emph{dot action} of the Weyl group $W$ of $\gfrak$ by setting
	$
		w\cdot\lambda\defeq w(\lambda+\rho)-\rho
	$
	for all $w\in W$ and $\lambda\in\hfrak^*$.  
	
\subsection{The Lie Algebras $\sllie_\infty$, $\solie_\infty$, and $\splie_\infty$}
\label{ch:glinfty}

The Lie algebras $\sllie_\infty$, $\solie_\infty$, and $\splie_\infty$ are, respectively, the direct limits of the finite-dimensional Lie algebras $\sllie_n$, $\solie_n$, and $\splie_{2n}$ with respect to root inclusions that increase the rank by $1$ (due to \cite{DP1999}, the Lie algebras do not depend on the choice of root inclusions).  That is, $\sllie_\infty=\lim\limits_{\underset{n}{\boldsymbol{\longrightarrow}}}\,\sllie_{n+1}$, $\solie_\infty=\lim\limits_{\underset{n}{\boldsymbol{\longrightarrow}}}\,\solie_{2n+1}=\lim\limits_{\underset{n}{\boldsymbol{\longrightarrow}}}\,\solie_{2n}$, and $\splie_\infty=\lim \limits_{\underset{n}{\boldsymbol{\longrightarrow}}}\,\splie_{2n}$. The remaining part of this subsection is based on \cite{BorelDC} and \cite{DPS}.

	Up to automorphism of $\sllie_\infty$, the Lie algebra $\sllie_\infty$ has a unique splitting maximal toral subalgebra $\hfrak_\textup{A}$.  There are two Dynkin Borel subalgebras $\bfrak$ containing $\hfrak_\textup{A}$ (up to automorphism of $\sllie_\infty$): $\bfrak_\textup{1A}$ and $\bfrak_\textup{2A}$.  The Borel subalgebra $\bfrak_\textup{1A}$ has the following Dynkin diagram $ \textup{A}_\infty^{\textup{1-sided}}$:
	\begin{align}
 \begin{dynkin}
    \dynkinline{1}{0}{2}{0};
    \dynkinline{2}{0}{3}{0};
    \dynkinline{3}{0}{4}{0};
    \dynkinline{4}{0}{5}{0};
    \dynkinline{5}{0}{6}{0};
    \dynkindots{6}{0}{7}{0};
    \foreach \x in {1,...,7}
    {
      \ifnum \x=7 {\fill (\dynkinstep*7,\dynkinstep*0) circle (0.01);}
       \else {
       \dynkindot{\x}{0}
       }
       \fi
    }
  \end{dynkin}\,,
\end{align}
while the Borel subalgebra $\bfrak_{\textup{2A}}$ is the Dynkin diagram $
  \textup{A}^{\textup{2-sided}}_\infty$: 
	\begin{align}\begin{dynkin}
    \dynkindots{1}{0}{2}{0};
    \dynkinline{2}{0}{3}{0};
    \dynkinline{3}{0}{4}{0};
    \dynkinline{4}{0}{5}{0};
    \dynkinline{5}{0}{6}{0};
    \dynkindots{6}{0}{7}{0};
    \foreach \x in {1,...,7}
    {
      \ifnum \ifnum \x=1 1 \else \ifnum \x=7 1 \else 0 \fi \fi =1 {\fill (\dynkinstep*\x,\dynkinstep*0) circle (0.01);}
       \else {
       \dynkindot{\x}{0}
       }
       \fi
    }
  \end{dynkin}\,.
\end{align}

	
For $\solie_\infty$, there are,up to automorphism of $\solie_\infty$, two splitting maximal toral subalgebras $\hfrak_\textup{B}$ and $\hfrak_\textup{D}$.  The splitting maximal toral subalgebra $\hfrak_\textup{B}$ corresponds to the direct limit $\solie_\infty = \lim\limits_{\underset{n}{\boldsymbol{\longrightarrow}}}\,\solie_{2n+1}$, whereas $\hfrak_\textup{D}$ corresponds to the direct limit $\solie_\infty =\lim\limits_{\underset{n}{\boldsymbol{\longrightarrow}}}\,\solie_{2n}$.  There is one (up to automorphism of $\solie_\infty$) Dynkin Borel subalgebra $\bfrak_\textup{B}$ containing $\hfrak_\textup{B}$ and has the following Dynkin diagram $\textup{B}_\infty$:
	\begin{align}
  \begin{dynkin}
    \dynkindoubleline{2}{0}{1}{0};
    \dynkinline{2}{0}{3}{0};
    \dynkinline{3}{0}{4}{0};
    \dynkinline{4}{0}{5}{0};
    \dynkinline{5}{0}{6}{0};
    \dynkindots{6}{0}{7}{0};
    \foreach \x in {1,...,7}
    {
      \ifnum \x=7 {\fill (\dynkinstep*7,\dynkinstep*0) circle (0.01);}
       \else {
       \dynkindot{\x}{0}
       }
       \fi
    }
  \end{dynkin}\,.
\end{align}
Up to automorphism of $\solie_\infty$, there is also one Dynkin Borel subalgebra $\bfrak_\textup{D}$ containing $\hfrak_\textup{D}$.  It has the following Dynkin diagram  $\textup{D}_\infty$:
	\begin{align}
 \begin{dynkin}
    \dynkinline{1}{0.5}{2}{0};
    \dynkinline{1}{-0.5}{2}{0};
    \dynkinline{2}{0}{3}{0};
    \dynkinline{3}{0}{4}{0};
    \dynkinline{4}{0}{5}{0};
    \dynkinline{5}{0}{6}{0};
    \dynkindots{6}{0}{7}{0};
    \foreach \x in {2,...,7}
    {
        \ifnum \x=7 {\fill (\dynkinstep*7,\dynkinstep*0) circle (0.01);}
       \else {
       \dynkindot{\x}{0}
       }
       \fi
    }
    {\fill (\dynkinstep*1,\dynkinstep*-0.5) circle (\dynkinradius);}
     {\fill (\dynkinstep*1,\dynkinstep*0.5) circle (\dynkinradius);}
  \end{dynkin}\,.
\end{align}

For $\splie_\infty$, there is only one splitting maximal toral subalgebra $\hfrak_\textup{C}$ (up to automorphism of $\splie_\infty$), and there is one (up to automorphism of $\splie_\infty$) Dynkin Borel subalgebra $\bfrak_\textup{C}$ containing $\hfrak_\textup{C}$ .  This Borel subalgebra has the following Dynkin diagram  $\textup{C}_\infty$:
	\begin{align}
  \begin{dynkin}
    \dynkindoubleline{1}{0}{2}{0};
    \dynkinline{2}{0}{3}{0};
    \dynkinline{3}{0}{4}{0};
    \dynkinline{4}{0}{5}{0};
    \dynkinline{5}{0}{6}{0};
    \dynkindots{6}{0}{7}{0};
    \foreach \x in {1,...,7}
    {
      \ifnum \x=7 {\fill (\dynkinstep*7,\dynkinstep*0) circle (0.01);}
       \else {
       \dynkindot{\x}{0}
       }
       \fi
    }
  \end{dynkin}\,.
  \end{align}
  

\subsection{Verma Modules}
\label{sec:verma}

	From now on, we assume that $\bfrak$ is a Dynkin Borel subalgebra of a root-reductive Lie algebra $\gfrak$ containing a splitting maximal toral subalgebra $\hfrak$.  Terms such as \emph{weight modules}, \emph{cyclic modules} and \emph{highest-weight modules} (with respect to $\bfrak$) are defined in the trivial manner.
	
	\begin{define}
	For every $\lambda\in\mathfrak{h}^*$, we define the \emph{Verma module} over $\mathfrak{g}$ of highest weight $\lambda$ to be the left $\Ulie(\mathfrak{g})$-module
		$
			\Mlie(\lambda;\mathfrak{g},\mathfrak{b},\mathfrak{h})\defeq \Ulie(\mathfrak{g})/I$,
		where $I$ is the left $\Ulie(\mathfrak{g})$-ideal generated by $\mathfrak{n}$ and $h-\lambda(h)$, for all $h \in \mathfrak{h}$.  If there is no confusion, we shall write $\Mlie(\lambda)$ for $\Mlie(\lambda;\mathfrak{g},\mathfrak{b,\mathfrak{h}})$. 
	\end{define}

	Note that $\Mlie(\lambda)$ is isomorphic to $\Ulie\left(\nfrak^-\right)$ as a left $\Ulie\left(\nfrak^-\right)$-module.  
	In addition, \cite{penkov-dimitrov} shows that $\Mlie(\lambda;\mathfrak{g},\mathfrak{b},\mathfrak{h})$ has a unique maximal proper $\Ulie(\mathfrak{g})$-submodule $N$, which is the sum of all proper submodules of $\Mlie(\lambda;\mathfrak{g},\mathfrak{b},\mathfrak{h})$.  Consequently,
	$
		\Llie(\lambda;\mathfrak{g},\mathfrak{b},\mathfrak{h})\defeq\Mlie(\lambda;\mathfrak{g},\mathfrak{b},\mathfrak{h})/N$
	is the unique simple quotient of $\Mlie(\lambda;\mathfrak{g},\mathfrak{b,\mathfrak{h}})$ (also denoted by $\Llie(\lambda)$ if there is no confusion).
	

	\begin{thm}
		For each $\lambda\in\hfrak^*$, the Verma module $\Mlie(\lambda)$ has at most one simple submodule, and if this submodule exists, it is also a Verma module.
		\label{thm:simplesubmodules}
	
	\begin{pf}
	Let $M:=\Mlie(\lambda)$.  First, observe that every submodule of $M$ has a singular vector with respect to $\bfrak$.  Suppose that $N_1$ and $N_2$ are $\Ulie(\mathfrak{g})$-submodules of $M$ with singular vectors $v_1 \neq 0$ and $v_2 \neq 0$, respectively. 
	Note that $M$ is isomorphic to $\Ulie\left(\mathfrak{n}^-\right)$ as an $\Ulie\left(\mathfrak{n}^-\right)$-module.  We can identify $M$ with $\Ulie\left(\mathfrak{n}^-\right)$, making $v_1$ and $v_2$ elements of $\Ulie\left(\mathfrak{n}^-\right)$.  Ergo, $N_1$ and $N_2$ are left $\Ulie\left(\mathfrak{n}^-\right)$-ideals $\Ulie\left(\mathfrak{n}^-\right)\cdot v_1$ and $\Ulie\left(\mathfrak{n}^-\right)\cdot v_2$, respectively.

	Because the subalgebra $\mathfrak{n}^-$ is locally finite, there exists a finite-dimensional subalgebra $\underline{\mathfrak{n}}^-$ of $\mathfrak{n}^-$ that contains the elements of $\mathfrak{n}^-$ involved in the PBW polynomial expressions for $v_1$ and $v_2$.  Consequently, the universal enveloping algebra $\Ulie\left(\underline{\mathfrak{n}}^-\right)$ is a noetherian ring.  
		From \cite[Lemma 4.1]{bggo}, we conclude that $\Ulie\left(\underline{\mathfrak{n}}^-\right)\cdot v_1$ must intersect $\Ulie\left(\underline{\mathfrak{n}}^-\right)\cdot v_2$ nontrivially.  Thence, $N_1$ and $N_2$  intersect nontrivially as well.  We then conclude that $N_1=N_2$.  In other words, $M$ has at most one simple submodule, and due to the fact that $\nfrak^-$ acts freely on $M$,  this submodule (if exists) is a Verma module.
	\end{pf}
	\end{thm}

	The theorem below gives a generalized version to the finite-dimensional case (see \cite{bggo}).  Theorem~\ref{thm:vermainf} offers an infinite-dimensional version of Verma's Theorem, which gives a condition under which there exists an embedding of a Verma module into another Verma module.
	
	\begin{thm}
		Let $\lambda,\mu\in\hfrak^*$.  The vector space $\Hom_{\Ulie(\mathfrak{g})}\big(\Mlie(\lambda),\Mlie(\mu)\big)$ is $0$- or $1$-dimensional.  Furthermore, all nonzero elements of  $\Hom_{\Ulie(\mathfrak{g})}\big(\Mlie(\lambda),\Mlie(\mu)\big)$ are embeddings.  If a nonzero homomorphism exists, then $\lambda \preceq \mu$.
		
	\begin{pf}
	Suppose that $\phi_1,\phi_2:\Mlie(\lambda)\to\Mlie(\mu)$ are nonzero $\mathfrak{g}$-module homomorphisms, where $\lambda,\mu\in\hfrak^*$.  We shall prove that $\phi_2=\kappa \, \phi_1$ for some $\kappa \in \amsbb{K}$.  Let $v_\lambda \neq 0$ and $v_\mu \neq 0$ be highest-weight vectors of $\Mlie(\lambda)$ and $\Mlie(\mu)$, respectively.  Write $w_i:=\phi_i\left(v_\lambda\right)$ for $i \in \{1,2\}$.   We identify $\Mlie(\mu)$ as a $\Ulie\left(\mathfrak{n}^-\right)$-module which is isomorphic to $\Ulie\left(\mathfrak{n}^-\right)$ itself.  Ergo, $w_1$ and $w_2$ are now elements of $\Ulie\left(\mathfrak{n}^-\right)$.   By the local finiteness of $\mathfrak{g}$, there exists a finite-dimensional  subalgebra $\underline{\mathfrak{g}}$ with a Borel subalgebra $\underline{\mathfrak{b}}:=\mathfrak{b}\cap\underline{\mathfrak{g}}$ that contains a maximal toral subalgebra $\underline{\mathfrak{h}}:=\mathfrak{h}\cap\underline{\mathfrak{g}}$.  Then, $\underline{M}:=\Ulie\left(\underline{\mathfrak{g}}\right)\cdot v_\mu$ is a Verma module over $\underline{\mathfrak{g}}$.  Now, $\Ulie\left(\underline{\mathfrak{g}}\right)\cdot w_1$ and $\Ulie\left(\underline{\mathfrak{g}}\right)\cdot w_2$ are isomorphic Verma modules over $\underline{\mathfrak{g}}$, both of which are embedded into $\tilde{M}$.  Since, in the finite-dimensional case, the homomorphism space between two Verma modules is either trivial or one-dimensional.  Therefore, we must have $\Ulie\left(\underline{\mathfrak{g}}\right)\cdot w_1=\Ulie\left(\underline{\mathfrak{g}}\right)\cdot w_2$.  Consequently, $w_2=\kappa\, w_1$ for some nonzero $\kappa \in \amsbb{K}$.  This means $\phi_2=\kappa\, \phi_1$.  Hence, $\textup{Hom}_{\Ulie(\mathfrak{g})}\big(\Mlie(\lambda),\Mlie(\mu)\big)$ is of dimension $0$ or $1$ over $\amsbb{K}$. 
	
	To show that any nonzero homomorphism in  $\textup{Hom}_{\Ulie(\mathfrak{g})}\big(\Mlie(\lambda),\Mlie(\mu)\big)$ must be an embedding, let $\phi$ be such a map.  Via the identification of $\Mlie(\lambda)$ and $\Mlie(\mu)$ with $\Ulie\left(\nfrak^-\right)$ as left $\Ulie\left(\nfrak^-\right)$ modules, we can easily see that $\phi$ is the multiplication map $x\mapsto u\boldsymbol{\cdot}x$ for some $ u \in \Ulie\left(\nfrak^-\right)$ and for all $x \in \Ulie\left(\nfrak^-\right)$.  Because $\Ulie\left(\nfrak^-\right)$ lacks zero divisors, $\phi$ must be injective.
	\end{pf}
	\label{thm:embeddings}
	\end{thm}
	
	\begin{define}
		Let $\lambda\in\hfrak^*$.  The global half sum of positive roots is denoted by $\rho$.
			\begin{enumerate}
				\item[(a)] We say that $\lambda$ is \emph{integral} if $\lambda\left(h_\alpha\right) \in \amsbb{Z}$ for every $\alpha \in \Delta$.
				\item[(b)] We say that $\lambda$ is \emph{antidominant} if $(\lambda+\rho)\left(h_\alpha\right) \notin \amsbb{Z}_{> 0}$ for any $\alpha \in \Delta^+$.
				\item[(c)] We say that $\lambda$ is \emph{almost antidominant} if $(\lambda+\rho)\left(h_\alpha\right) \in \amsbb{Z}_{> 0}$ for only finitely many $\alpha \in \Delta^+$.
			\end{enumerate}
			\label{def:wttypes}
		\end{define}
	
	\begin{thm}[Verma's Theorem]
			For $\lambda \in \hfrak^*$ and for a given a positive root $\alpha$ such that $s_\alpha\cdot \lambda \preceq \lambda$, there exists an embedding $\Mlie\left(s_\alpha\cdot\lambda\right)\overset{\subseteq}{\longrightarrow}\Mlie(\lambda)$.
			\begin{pf}
				For $n \in \amsbb{Z}_{>0}$, write $\bfrak_n$ and $\hfrak_n$ for $\bfrak\cap\gfrak_n$ and $\hfrak\cap\gfrak_n$, respectively.   Let $\lambda_n$ be the restriction of $\lambda$ onto $\hfrak_n$.  Denote by $M$ the Verma module $\Mlie(\lambda;\gfrak,\bfrak,\hfrak)$, while $M_n$ is the Verma module $\Mlie\left(\lambda_n;\gfrak_n,\bfrak_n,\hfrak_n\right)$.  If $u$ is a highest-weight vector of $M$, then by identifying a highest-weight vector of $M_n$ with $u$, we have $M_1 \subseteq M_2 \subseteq M_3 \subseteq \ldots$.  Clearly, as a set, $M$ is the direct limit of $\left(M_n\right)_{n\in\amsbb{Z}_{>0}}$ under inclusion maps.
				
				The root space $\gfrak^{-\beta}$ is finite dimensional for every $\beta \in \Psi:=\left\{\gamma\in\Delta^+\suchthat{} \gamma \preceq \lambda -s_\alpha\cdot\lambda\right\}$.  The set $\Psi$ is clearly finite.  Therefore, for sufficiently large values of $n$, say $n \geq m$ for some $m \in \amsbb{Z}_{>0}$, we have $\gfrak^{-\beta}\subseteq \gfrak_n$ for all $\beta \in \Psi$, which further implies that $\alpha|_{\hfrak_n}$ is a positive root of $\gfrak_n$, and $s_{\alpha} \in W_n$.  Thus, for such $n\in\amsbb{Z}_{>0}$, the Verma module $\Mlie\left(s_\alpha\cdot\lambda_n;\gfrak_n,\bfrak_n,\hfrak_n\right)$ is isomorphic to a unique $\gfrak_n$-submodule $N_n$ of $M_n$, where we have applied the finite-dimensional version of Verma's Theorem.   

				Now, observe that, for $n \geq m$, $s_\alpha \cdot \lambda_n \in \hfrak_n^*$ is identical to the restriction of $s_\alpha \cdot\lambda_{n+1} \in \hfrak_{n+1}^*$ onto $\hfrak_n$.  Furthermore, the weight space associated to the weight $s_\alpha\cdot\lambda_n$ of $M_n$ (where the dot action is done in $\hfrak_n^*$) is precisely $M^{s_\alpha\cdot\lambda}$.    This means that the highest-weight spaces of $N_n$ and of $N_{n+1}$, which correspond to the weights $s_\alpha\cdot\lambda_n$ and $s_\alpha\cdot\lambda_{n+1}$, respectively, are identical for $n \geq m$.  That is, $N_n \subseteq N_{n+1}$ for every integer $n \geq m$.  The direct limit $N:=\lim_{\underset{n}{\boldsymbol{\longrightarrow}}} N_n$ of $\left(N_n\right)_{n  \geq m}$ under inclusion maps is thus a $\gfrak$-submodule of $M$ isomorphic to $\Mlie(s_\alpha\cdot\lambda;\gfrak,\bfrak,\hfrak)$.
			\end{pf}
			
			\label{thm:vermainf}
	\end{thm}

	It turns out that the BGG Theorem for the finite-dimensional case have infinite-dimensional analogues.   The generalizations below shall be called the BGG Theorem as well.
	
	\begin{thm}[BGG Theorem]
		For $\lambda,\mu\in\hfrak^*$, there exists a nontrivial $\gfrak$-module homomorphism from $\Mlie(\lambda)$ to $\Mlie(\mu)$ if and only if $\lambda$ is \emph{strongly linked} to $\mu$, namely, there exist positive roots $\alpha_1,\alpha_2,\ldots,\alpha_l$ with $l$ being a nonnegative integer such that
		\begin{align}
			\lambda = \left(s_{\alpha_l}\cdots s_{\alpha_2}s_{\alpha_1}\right)\cdot\mu \preceq \left(s_{\alpha_{l-1}}\cdots s_{\alpha_2}s_{\alpha_1}\right)\cdot\mu \preceq \ldots \preceq s_{\alpha_1}\cdot\mu \preceq \mu\,.
			\label{eq:stronglinkage}
		\end{align}
		That is, for $\mu \in \hfrak^*$, all Verma submodules of $\Mlie(\mu)$ is of the form $\Mlie(w\cdot \mu)$, where $w$ is an element of the Weyl group.
		\begin{pf}
			The converse is clear, so we prove the direct implication.  Suppose that $\Mlie(\lambda)$ is a submodule $N$ of $M:=\Mlie(\mu)$.  Let $u$ and $v$ be highest-weight vectors of $M$ and $N$, respectively.  For each $n \in \amsbb{Z}_{>0}$, write $M_n$ for $\Ulie\left(\gfrak_n\right)\cdot u$.  As $M=\lim_{\underset{n}{\boldsymbol{\longrightarrow}}}M_n$, there exists $m \in \amsbb{Z}_{>0}$ such that $n \geq m$ implies $v \in M_n$.  
			
			For $n \geq m$, write $N_n$ for $\Ulie\left(\gfrak_n\right)\cdot v$.  Then, $N_n$ is a Verma submodule of $M_n$ (over $\gfrak_n$).  The finite-dimensional BGG Theorem guarantees that $\lambda_n:=\lambda|_{\hfrak_n}$ is strongly linked to $\mu_n:=\mu|_{\hfrak_n}$.  The positive roots $\alpha_n^j$, $j=1,2,\ldots,l_n$, involved in the $n$-th linkage
			\begin{align}
				\lambda_n  = \left(s_{\alpha_n^{l_n}}\cdots s_{\alpha_n^2}s_{\alpha_n^1}\right)\cdot \mu_n \preceq \left(s_{\alpha_n^{l_n-1}}\cdots s_{\alpha_n^2}s_{\alpha_n^1}\right)\cdot \mu_n \preceq \ldots \preceq s_{\alpha_n^1}\cdot \mu_n \preceq \mu_n
			\end{align} 
			between $\lambda_n$ and $\mu_n$ must belong to the set $\left\{\alpha \in \Delta^+\suchthat{} \alpha \preceq \mu-\lambda\right\}$, which is a finite set.  
			
			If $\mu-\lambda =\sum_{\alpha \in \Sigma}t_\alpha \alpha$ and $t_\alpha \in \amsbb{Z}_{\geq 0}$ for each $\alpha \in \Sigma$, then the lenth $l_n$ of the $n$-th linkage is at most $\sum_{\alpha \in \Sigma}t_\alpha <\infty$.  Using the Pigeonhole Principle, it follows that there are infinitely many $n \geq m$ with the same linkage pattern, say
			$
				\lambda_n  = \left(s_{\alpha_{l}}\cdots s_{\alpha_2}s_{\alpha_1}\right)\cdot \mu_n \preceq \left(s_{\alpha_{l-1}}\cdots s_{\alpha_2}s_{\alpha_1}\right)\cdot \mu_n \preceq \ldots \preceq s_{\alpha_1}\cdot \mu_n \preceq \mu_n$,
			where $\alpha_1,\alpha_2,\ldots,\alpha_l$ are positive roots.  Hence, Equation (\ref{eq:stronglinkage}) holds, and the proof is complete.
		\end{pf}
		\label{thm:stronglinkagegen}
	\end{thm}

	\begin{thm}
			For $\lambda \in\hfrak^*$, $\Mlie(\lambda)$ is simple if and only if $\lambda$ is antidominant.
			
			\begin{pf}
				Let $\rho$ be a global half sum of positive roots.  For each root $\alpha$, $s_\alpha$ is the reflection with respect to $\alpha$ and $h_{\alpha}$ is as defined in Definition~\ref{def:wttypes}.
				
				\begin{itemize}
					\item[($\Rightarrow$)]  Suppose that $\lambda\in\hfrak^*$ is not an antidominant weight. Then, there exists a positive root $\alpha$ such that $(\lambda+\rho)\left(h_\alpha\right) \in \amsbb{Z}_{>0}$.  This means $s_\alpha \cdot \lambda \precneq \lambda$ and
					$
						0 \subsetneq \Mlie\left(s_\alpha\cdot\lambda\right)\subsetneq \Mlie(\lambda)
					$, where we have applied Verma's Theorem (\ref{thm:vermainf}); as a result, $\Mlie(\lambda)$ is not simple.
					\item[($\Leftarrow$)]  
					Suppose that $\Mlie(\lambda)$ is not simple.  Then, it has a proper nonzero submodule, which must have a highest-weight vector whose weight is $\mu \in \hfrak^*$.  Then, $\Mlie(\mu)$ is a proper Verma submodule of $\Mlie(\lambda)$, so $\mu \precneq \lambda$.  Using the BGG Theorem (\ref{thm:stronglinkagegen}), there are positive roots $\alpha_1$, $\alpha_2$, $\ldots$, $\alpha_l$ with $l \in \amsbb{Z}_{>0}$ such that $
						\mu= \left(s_{\alpha_l}\cdots s_{\alpha_2}s_{\alpha_1}\right)\cdot\lambda \preceq \left(s_{\alpha_{l-1}}\cdots s_{\alpha_2}s_{\alpha_1}\right)\cdot\lambda \preceq \ldots \preceq s_{\alpha_1}\cdot\lambda \preceq \lambda$.
					Because $s_{\alpha_1}\cdot\lambda \preceq \lambda$, we have $(\lambda+\rho)\left(h_{\alpha_1}\right) \in \amsbb{Z}_{>0}$.  Thence, $\lambda$ is not antidominant.
				\end{itemize}
			\end{pf}
			\label{thm:simpleverma}
		\end{thm}
		
		\begin{thm}
			Let $\lambda  \in \hfrak^*$.  The following conditions are equivalent:
			\begin{itemize}
				\item[(a)]  The module $\Mlie(\lambda)$ is of finite length.
				\item[(b)] The module $\Mlie(\lambda)$ has a simple submodule.
				\item[(c)] There exists an antidominant weight $\xi$ such that $\xi$ is strongly linked to $\lambda$.
				\item[(d)] The weight $\lambda$ is almost antidominant.
			\end{itemize}
			
			\begin{pf}
				For simplicity, we shall denote $M$ for $\Mlie(\lambda)$.
				
				\begin{itemize}
				\itemindent1cm
					\item[$\big($(a)$\Leftrightarrow$(b)$\big)$]  For the direct implication, let 
						$
							0=M_0\subsetneq M_1 \subsetneq M_2 \subsetneq \ldots \subsetneq M_{l-1}\subsetneq M_l=M
						$ be a composition series of $M$ for some $l \in \amsbb{Z}_{>0}$.  Then, $M_1$ must be simple.  
						
						Conversely, let $L$ be a simple submodule of $M$.  Then, $L$ is a Verma module of highest weight $\mu \in \hfrak^*$ with $\mu \preceq \lambda$.  Note that every nonzero submodule of $M$ must include $L$.  Any composition series $0=M_0\subsetneq M_1 \subsetneq M_2  \subsetneq  \ldots \subsetneq M_{l-1}\subsetneq M_l=M$ by submodules of $M$ must have $L \subseteq M_1 \subsetneq M_2 \subsetneq \ldots \subsetneq M_l$.  We look at the set $S_i$ of all possible highest-weight spaces of $M_i$, for $i=1,2,\ldots,l$.  Clearly, for every $i=1,2,\ldots,l$, $\mu \preceq \xi$, where $\xi$ is the weight of a highest-weight space in $S_i$.  There are only finitely many weights $\xi$ for which $\mu \preceq \xi \preceq \lambda$, and  the weight space with weight $\xi$ is finite dimensional for each $\xi \in \hfrak^*$.  If $m$ is the sum of the dimensions of all the weight spaces with weight $\xi \in \hfrak^*$ such that $\mu \preceq \xi \preceq \lambda$, then we have that $l \leq m$.  Hence, $M$ is of finite length.
					\item[$\big($(b)$\Leftrightarrow$(c)$\big)$] We can easily apply Verma's Theorem (Theorem \ref{thm:vermainf}), the BGG Theorem (Theorem \ref{thm:stronglinkagegen}), and Theorem~\ref{thm:simpleverma} to verify that (b) and (c) are equivalent. 
					
						\item[$\big($(a)$\Leftarrow$(d)$\big)$]
						Suppose $L$ is a simple submodule of $M$.  Then, $L$ is a Verma module with the highest weight $\mu \preceq \lambda$, for some $\mu \in \hfrak^*$.  By the BGG Theorem (Theorem \ref{thm:stronglinkagegen}), $w\cdot\lambda = \mu$ for some $w \in W$.  Since $\mu$ is antidominant (by Theorem~\ref{thm:simpleverma}), this means $\lambda$ is almost antidominant.
						
						\item[$\big($(d)$\Leftarrow$(c)$\big)$]	Suppose that $\lambda$ is almost antidominant.  For each $\mu\in\hfrak^*$, let $\Omega(\mu)$ denotes the set of positive roots $\alpha$ such that $h_\alpha(\lambda+\rho)$ is a positive integer.  We say that $\alpha \in \Omega(\mu)$ is \emph{minimal} if $\alpha$ cannot be written as a sum of at least two elements of $\Omega(\mu)$.  Let $\omega(\mu)$ be the cardinality of $\Omega(\mu)$.  
		 	
		 	By the assumption, $\omega(\lambda)<\infty$.  We shall prove by induction on $\omega(\lambda)$.  Pick a minimal $\alpha\in \Omega(\lambda)$.  Then, we have $s_\alpha\cdot \lambda\preceq \lambda$ and $\omega(s_\alpha\cdot \lambda)<\omega(\lambda)$.  By the induction hypothesis, there exists an antidominant weight $\xi$ such that $\xi$ is strongly linked to $s_\alpha\cdot\lambda$.  That is, there are positive roots $\alpha_2,\alpha_3,\ldots,\alpha_l$ such that 
				\begin{align}
					\xi = \left(s_{\alpha_l}\cdots s_{\alpha_2}\right)\cdot\left(s_{\alpha_1}\cdot \lambda\right) \preceq \left(s_{\alpha_{l-1}}\cdots s_{\alpha_2}\right)\cdot\left(s_{\alpha_1}\cdot \lambda\right) \preceq \ldots \preceq s_{\alpha_2}\cdot\left(s_{\alpha_1}\cdot \lambda\right) \preceq s_{\alpha_1}\cdot\lambda
				\end{align}
				It follows immediately that				
				$
			\xi = \left(s_{\alpha_l}\cdots s_{\alpha_2}s_{\alpha_1}\right)\cdot\lambda \preceq \left(s_{\alpha_{l-1}}\cdots s_{\alpha_2}s_{\alpha_1}\right)\cdot\lambda \preceq \ldots \preceq s_{\alpha_1}\cdot\lambda \preceq \lambda$,
			and our proof is now complete.
				\end{itemize}
			\end{pf}
			\label{thm:finitelength}
		\end{thm}

\section{The Extended Categories $\bggO$}
\label{ch:bbggo}

In this section, $\gfrak$ is a root-reductive algebra.  The subalgebra $\bfrak$ is a Dynkin Borel subalgebra of $\gfrak$ containing a splitting maximal toral subalgebra $\hfrak$.

\subsection{The Definition}
\label{sec:extbggo}

\begin{define}
	The \emph{extended category $\bggO$}, denoted by $\bar\bggO^\gfrak_\bfrak$, is the full subcategory of the category of $\gfrak$-modules satisfying the following two conditions:
	\begin{enumerate}[(i)]
		\item  Every $M\in\bar\bggO^\gfrak_\bfrak$ is an $\hfrak$-weight $\gfrak$-module with finite-dimensional $\hfrak$-weight spaces;
		\item Every $M\in\bar\bggO^\gfrak_\bfrak$ is locally $\nfrak$-finite (that is, $\Ulie(\nfrak)\cdot v$ is finite dimensional for every $v\in M$).
	\end{enumerate}
\end{define}

When this cannot cause confusion, we shall write $\bar\bggO$ for $\bar\bggO^{\gfrak}_\bfrak$.  We can define the duality functor of the category $\bbggO$ in the same manner as the standard duality of the category $\bggO$.  More precisely, we have the following definition.

\begin{define}
	Let $M\in\bar\bggO$.  For $\lambda\in\hfrak^*$, let $M^{\vee,\lambda}$ denote the set of $f\in M^*$ such that $f$ vanishes on $M^\mu$ for every $\mu\in\hfrak^*\smallsetminus\{\lambda\}$.  The \emph{dual} of $M$ is defined to be $M^\vee\defeq \bigoplus_{\lambda\in\hfrak^*}\,M^{\vee,\lambda}$.
\end{define}

Now, if $\left\{x_{\pm\alpha}\,\boldsymbol{|}\,\alpha\in\Delta^+\right\}\cup \left\{h_\beta\,\boldsymbol{|}\,\beta\in\Sigma^+\right\}$ is a \emph{Chevalley basis} (see \cite{humphreys} and \cite{bggo}) of $\gfrak$, then the \emph{transpose map} $\tau:\gfrak\to\gfrak$ is the linear map given by $\tau(h):=h$ for all $h\in\hfrak$, and $\tau\left(x_{\pm\alpha}\right):=x_{\mp\alpha}$ for all positive roots $\alpha$.  Note that $\big[\tau(x),\tau(y)\big]=-\tau\big([x,y]\big)$ for all $x,y\in\gfrak$.  We have the following proposition, which can be proven in the same way as in the finite-dimensional setting.

\begin{prop}
	For every object $M\in\bar\bggO$, $M^\vee$ is a $\gfrak$-module with respect to the twisted $\gfrak$-action $(g\cdot f)(v)\defeq f\big(\tau(g)\cdot v\big)$ for all $g\in\gfrak$, $v\in M$, and $f\in M^\vee$.  Furthermore, $M^\vee \in \bar\bggO$.
\end{prop}
	
	Finally, we consider the categories $\bggO$ and $\bbggO$ for a finite-dimensional reductive Lie algebra $\gfrak$. With $\bar\bggO$ and $\bggO$ being highest-weight categories, we automatically have BGG reciprocity (which is a special case of Brauer-Humphreys reciprocity \cite[Theorem 3.11]{hwcat}):
	\begin{align}
		\big[\Vlie(\mu):\Llie(\lambda)\big]=\big[\Mlie(\mu):\Llie(\lambda)\big]=\big\{\Ilie(\lambda):\Vlie(\mu)\big\}=\big\{\Plie(\lambda):\Mlie(\mu)\big\}\,,
	\end{align}
	where $\Vlie(\mu)$ is the dual Verma module $\big(\Mlie(\mu)\big)^\vee$, $\Ilie(\lambda)$ is the injective hull of the simple object $\Llie(\lambda)$, and $\Plie(\lambda)=\big(\Ilie(\lambda)\big)^\vee$ is the projective cover of $\Llie(\lambda)$.  However, as we shall later prove, the category $\bar\bggO^\gfrak_\bfrak$ is not a highest-weight category if $\gfrak$ is infinite dimensional.

\subsection{Direct Sum Decompositions}

Let $\gfrak$ be a root-reductive Lie algebra.  The objective of this section is to verify that every object in $\bar\bggO$ has a decomposition into a direct sum of indecomposable objects.  Furthermore, this decomposition is unique up to isomorphism.  That is, if an object $M\in\bar\bggO$ can be written as $\bigoplus_{j\in J}\,M_j=M=\bigoplus_{j'\in J'}\,M_{j'}'$, where $J$ and $J'$ are index sets and $M_j,M'_{j'}\in\bar\bggO$ are indecomposable for all $j\in J$ and $j'\in J'$, then there exists a bijection $\psi:J\to J'$ such that $M_j\cong M'_{\psi(j)}$ for every $j\in J$.

First, we need the proposition below.  The proof is trivial.

\begin{prop}
	Every indecomposable object $M\in\bar\bggO$ satisfies $\textup{supp}(M)\subseteq \lambda+\Lambda$, where $\Lambda$ is the root lattice.  In particular, $\textup{supp}(M)$ is countable and $M$ is countable dimensional.
	\label{prop:indecobj_weight}
\end{prop}

\begin{thm}
	Every $M\in\bar\bggO$ is a direct sum of indecomposable objects.
	\label{thm:directsumdecomposition}
	
	\begin{pf}
		For an object $M\in\bar\bggO$, we say that $\xi\in\hfrak^*$ is a \emph{decomposable weight} of $M$ if there exist submodules $N$ and $L$ of $M$ such that $M=N\oplus L$ with $\dim_\amsbb{K}\left(N^\xi\right)>0$ and $\dim_\amsbb{K}\left(L^\xi\right)>0$.  Otherwise, $\xi$ is an \emph{indecomposable weight}.  (Note that, by abuse of language, an indecomposable weight of a $\gfrak$-module $X$ need not be a weight of $X$.  In other words, if $X^\xi=0$, then $\xi$ is an indecomposable weight of $X$, despite not actually being a weight of $X$.)

		For a semisimple $\hfrak$-module $X$, the \emph{support} $\textup{supp}(X)$ of $X$ is the set of the $\hfrak$-weights of $X$.  For a subset $S\subseteq\hfrak^*$, we say that $S$ is an \emph{indecomposable weight set} of $M$ if every weight in $S$ is an indecomposable weight of $M$ and if $M$ cannot be written as a direct sum $M=N\oplus L$ such that $\textup{supp}(N)\cap S$ and $\textup{supp}(L)\cap S$ are both nonempty.
	
		Assume that $M\in\bar\bggO$ is nonzero.   From Proposition~\ref{prop:indecobj_weight} above, we may assume that $M$ is countable dimensional.  Hence, $\textup{supp}(M)=\left\{\mu_1,\mu_2,\mu_3,\ldots\right\}$ for some weights $\mu_1,\mu_2,\mu_3,\ldots \in \hfrak^*$. 
		
		We shall prove by (countable) transfinite induction that, for each $i=1,2,3,\ldots$, there exists an index set $J_i$ such that $M$ has a direct sum decomposition
		$
			M=\bigoplus_{j\in J_i}\,D_i(j)$
		such that the set $Q_i:=\big\{\mu_1,\mu_2,\ldots,\mu_{i}\big\}$ is an indecomposable weight set of $D_i(j)$ (here, $Q_1:=\emptyset$).  We further require that the decomposition above (with respect to the weight $\mu_i$) be compatible with the decomposition with respect to every weight $\mu_{i'}\in Q_i$ in the sense that, for any $j\in J_i$, there exists (uniquely) $j'\in J_{i'}$ such that
		$
			D_i(j)\subseteq D_{i'}\left(j'\right)$.
		
		For the base case, if $\mu_1$ is already an indecomposable weight, then $M=M$ is a required decomposition (that is, $J_1=\{1\}$ and $D_1(1)=M$).  If $\mu_1$ is a decomposable weight, then there exists a submodule $D_1(1)$ of $M$ such that $D_1(1)$ is a direct summand of $M$.  We may chose $D_1(1)$ so that $\dim_\amsbb{K}\Big(\big(D_1(1)\big)^{\mu_1}\Big)>0$ is minimal.  Then, $\mu_1$ is an indecomposable weight of $D_1(1)$.  Let $D'_1(1)$ denote a complementary submodule of $M$ to $D_1(1)$.  We proceed further by decomposing $D'_1(1)$ as a direct sum of submodules.  As $\dim_\amsbb{K}\left(M^{\mu_1}\right)<\infty$, the process will lead after finitely many steps to a decomposition
		$
			M=D_1(1)\oplus D_1(2)\oplus\ldots\oplus D_1\left(n\right)
		$
		such that $\mu_1$ is an indecomposable weight of each $D_1(j)$ for $j=1,2,\ldots,n$.  Thus, $Q_1=\left\{\mu_1\right\}$ is an indecomposable weight set of each $D_1(j)$, and we set $J_1:=\{1,2,\ldots,n\}$.
		
		Suppose that the induction hypothesis holds for $\mu_1,\mu_2,\ldots,\mu_{i-1}$.  That is, we have a direct sum decomposition
		$
			M=\bigoplus_{j\in J_{i'}}\,D_{i'}(j)$
		such that $Q_{i'}$ is an indecomposable weight set of $D_{i'}(j)$, where $i'=1,2,\ldots,i-1$.  We proceed to decompose $D_{i-1}(j)$ with respect to the weight $\mu_i$ instead of $\mu_1$ in the same manner as the base case.  That is, 
		$
			D_{i-1}(j)=\bigoplus_{k\in J_{i-1}^j}\,D_{i-1}^j\left(k\right)$
		for some index set $J_{i-1}^j$ and for some submodules $D_{i-1}^j(k)$ such that $\mu_i$ is an indecomposable weight.  Therefore,
		$
			M=\bigoplus_{j\in J_{i-1}}\,\bigoplus_{k\in J_{i-1}^j}\,D_{i-1}^j(k)
		$
		is a decomposition in which $Q_{i-1}$ is an indecomposable weight set of each $D_{i-1}^j(k)$ and $\mu_i$ is an indecomposable weight of each $D_{i-1}^j(k)$.  
		
		Note that, for a given submodule $D_{i-1}^j(k)$, $Q_i$ may not be an indecomposable weight set of $D_{i-1}^j(k)$, but when that is the case, we can further decompose $D_{i-1}^j(k)$ as follows:
		\begin{align}
			D_{i-1}^j(k)= \bar{D}_{i-1}^j(k)\oplus \tilde{D}_{i-1}^j(k)\,,
		\end{align}
		where $\left(\bar{D}_{i-1}^j(k)\right)^{\mu_i}=0$ and $\left(\tilde{D}_{i-1}^j(k)\right)^{\mu_{i'}}=0$ for all $i'=1,2,\ldots,i-1$.  Let $\tilde{J}_{i-1}^j$ denote the subset of $J^j_{i-1}$ consisting of $k\in J_{i-1}^j$ such that $Q_i$ is not an indecomposable weight set of $D_{i-1}^j(k)$.  Then,
		\begin{align}
			M=\bigoplus_{j\in J_{i-1}}\,\left(\left(\bigoplus_{k\in J_{i-1}^j\smallsetminus\tilde{J}_{i-1}^j}\,D_{i-1}^j(k) \right)\oplus\left(\bigoplus_{k\in \tilde{J}_{i-1}^j}\,\left(\bar{D}_{i-1}^j(k)\oplus\tilde{D}_{i-1}^j(k)\right)\right)_{\vphantom{1_1}}^{\vphantom{1^1}}\right)
			\label{eq:hugedirectsum}
		\end{align}
		is a direct sum decomposition of $M$ with respect to weight $\mu_i$ and with the required properties.
		
		To complete the proof, we define the partially ordered set $\mathscr{P}$ to be the set of all pairs $\big(i,D_i(j)\big)$ where $\mu\in\textup{supp}(M)$ and $j\in J_\mu$, equipped with the partial order $\preceq$ defined by 
		\begin{align}
		\Big(i,D_i(j)\Big) \preceq \Big(i',D_{i'}\left(j'\right)\Big)\textup{ if and only if }i\leq i'\textup{ and }D_i(j)\supseteq D_{i'}\left(j'\right)\,.
		\end{align}  
		We write $\mathscr{M}$ for the set of maximal chains in $\mathscr{P}$.   Using the finite-dimensionality assumption on the weight spaces of $M$, we conclude that
		\begin{align}
			M=\bigoplus_{\mathscr{C}\,\in\,\mathscr{M}}\,\mathcal{D}(\mathcal{C})\,,\textup{ where }
			\mathcal{D}(\mathcal{C}):=\bigcap_{\big(i,D_i\left(j_i\right)\big)\in\mathcal{C}}\,D_i\left(j_i\right)
		\end{align} 
		for every maximal chain $\mathcal{C}$ in $\mathscr{P}$.  For each $i=1,2,3,\ldots$, the set $Q_i$ is an indecomposable weight set of each $D_i\left(j_i\right) \in\mathcal{C}$, for a given maximal chain $\mathcal{C}\subseteq \mathscr{P}$.  Then, $\textup{supp}(M)=\bigcup_{i\in\mathbb{Z}_{>0}}\,Q_i$ is an indecomposable weight set of every $\mathcal{D}(\mathcal{C})$.  Consequently, $\mathcal{D}(\mathcal{C})$ is an indecomposable module.  (Note that a direct summand $\mathcal{D}(\mathcal{C})$ may be the zero module for some $\mathcal{C}\in\mathcal{M}$, but this does not affect the proof or the statement of this theorem.)
	\end{pf}
\end{thm}

\begin{prop}
	Let $M\in\bar\bggO$ be indecomposable.  Then, every $\varphi\in \End_{\bar\bggO}(M)$ is either an automorphism or is \emph{locally nilpotent} (namely, for every $v\in M$, there exists $k\in\amsbb{Z}_{\geq 0}$ such that $\varphi^{ k}(v)=0$).
	
	\begin{pf}
		Let $K_k:=\ker\left(\varphi^{k}\right)$ and $I_k:=\im\left(\varphi^{k}\right)$ for each $k=0,1,2,\ldots$ (here, $\varphi^0$ is the identity map $\textup{id}_M$).  We observe that $K_0\subseteq K_1\subseteq K_2 \subseteq \ldots$ and $I_0\supseteq I_1\supseteq I_2\supseteq \ldots$.  Set $K:=\bigcup_{k=0}^\infty\,K_k$ and $I:=\bigcap_{k=0}^\infty\,I_k$.  
		
		Fix $\lambda\in\hfrak^*$.  The restriction $\psi_\lambda$ of $\varphi$ onto $M^\lambda$ is a linear map on a finite-dimensional vector space.  Hence, $M^\lambda$ decomposes as $\im\left(\psi_\lambda^{ k}\right)\oplus \ker\left(\psi_\lambda^{k}\right)$ for every $k=0,1,2,\ldots$.  Since $M^\lambda$ is a finite-dimensional vector space and 
		$
			\im\left(\psi_\lambda\right)\supseteq \im\left(\psi_\lambda^{2}\right)\supseteq \im\left(\psi_\lambda^{3}\right)\supseteq \ldots$, the submodules $\im\left(\psi_\lambda^{k}\right)$, where $k=0,1,2,\ldots$, must stabilize.  Assume that, for some $j\in\amsbb{Z}_{\geq 0}$, we have 
	\begin{align}\im\left(\psi_\lambda^{ j}\right)=\im\left(\psi_\lambda^{j+1}\right)=\im\left(\psi_\lambda^{j+2}\right)=\ldots\,.
	\end{align}  That is, the kernels must also stabilize at the same index: 
	\begin{align}
	\ker\left(\psi_\lambda^{ j}\right)=\ker\left(\psi_\lambda^{j+1}\right)=\ker\left(\psi_\lambda^{j+2}\right)=\ldots\,.
	\end{align}  This shows that $K^\lambda$ and $I^\lambda$ are equal to $K_j^\lambda$ and $I_j^\lambda$ for some $j\in\amsbb{Z}_{\geq 0}$, depending on $\lambda$.  Therefore, the sum 
	$
		(K+L)^\lambda=K^\lambda+L^\lambda =K_j^\lambda+I_j^\lambda
	$ is direct and equals $M^\lambda$, as $M^\lambda =\im\left(\psi_\lambda^{ j}\right)\oplus\ker\left(\psi_\lambda^{j}\right)=I_j^\lambda\oplus K_j^\lambda$.  Since $\lambda$ is arbitrary, we obtain $M=K\oplus I$.
		
		As $M$ is an indecomposable object, we have either $K=0$ or $I=0$.  In the former case, we conclude that $\varphi$ is an isomorphism, and in the latter case, we see that $\varphi$ is locally nilpotent.
	\end{pf}
\end{prop}

\begin{prop}
	For every indecomposable object $M\in\bar\bggO$, the endomorphism ring $\End_{\bar\bggO}{(M)}$ is a local ring.
	
	\begin{pf}
		Let $J\subseteq R:=\End_{\bar\bggO}(M)$ be the set of all locally nilpotent endomorphisms of $M$.  By the previous proposition, $J$ is the set of all non-invertible elements of $R$.  We must prove that $J$ is an ideal of $R$.  
		
		First, if $\varphi\in J$ and $\psi\in R$, then $\varphi\circ\psi$ cannot be an epimorphism because $\varphi$ is not surjective on any weight space of $M$, and $\psi\circ\varphi$ is not a monomorphism because $\varphi$ is not injective on any weight space of $M$.  That is, $\varphi\circ\psi$ and $\psi\circ\varphi$ are both in $J$.  
		
		Now, we assume that $\varphi,\psi\in J$.  We must show that $\varphi+\psi$ belongs to $J$ too.  Suppose on the contrary that $\varphi+\psi\notin J$.  Then, $\varphi+\psi$ is invertible.  Hence, $\varphi+\psi=\phi$ for some automorphism $\phi:M\to M$.  Let $\alpha:=\varphi\circ\phi^{-1}$ and $\beta:=\psi\circ\phi^{-1}$.  Then,  $\alpha+\beta=\id_M$ and $\alpha,\beta\in J$.  Note that $\alpha\circ\beta=\beta\circ\alpha$.  Fix a weight $\lambda$ of $M$.  Suppose that $\alpha^{k}$ and $\beta^{l}$ vanish on $M^\lambda$, for some $k,l\in\amsbb{Z}_{> 0}$.  Then, $
			(\alpha+\beta)^{ k+l}=\sum_{r=0}^{k+l}\,\binom{k+l}{r}\,\alpha^{ r}\circ \beta^{k+l-r}$ must vanish on $M^\lambda$ as well.  Ergo, the endomorphism $\alpha+\beta$ cannot equal $\id_M$, which is a contradiction.
	\end{pf}
\end{prop}

From the proposition above, the Krull-Schmidt-Remak-Azumaya Theorem (see \cite{KSRA}) immediately implies the following corollary.

\begin{cor}
	Every object in $\bar\bggO$ admits a unique, up to isomorphism, decomposition into a direct sum of indecomposable objects.
	\label{cor:KSRA}
\end{cor}
	
\subsection{Generalized Composition Series}

In this subsection, we shall prove that every module in $\bbggO$ has an analogue of composition series.  Following the proof by V. Kac of~\cite[Lemma 9.6]{kac}, we have the following lemma.

\begin{lem}
	Let $M\in\bar\bggO$ and $\lambda\in\hfrak^*$.  Suppose that all weights $\xi$ of $M$ satisfy $\xi\preceq \upsilon$ for some fixed upper bound $\upsilon\in\hfrak^*$.  Then, there exist a $\gfrak$-module filtration 
	\begin{align}
		0=M_0\subseteq M_1\subseteq M_2\subseteq \ldots \subseteq M_{k-1}\subseteq M_k=M
	\end{align} 
	and a subset $J\subseteq \{1,2,\ldots,k\}$ such that
	\begin{enumerate}[(i)]
		\item if $j\in J$, then $M_{j}/M_{j-1}\cong \Llie\big(\xi(j)\big)$ for some $\xi(j)\in \hfrak^*$ with $\xi(j)\succeq \lambda$,
		\item if $j\notin J$, then $\left(M_{j}/M_{j-1}\right)^\mu=0$ for every $\mu\succeq \lambda$.
	\end{enumerate}
	\label{lem:kac_series}
	
%
\end{lem}


\begin{cor}
	Let $M\in\bar\bggO$ and $\lambda,\nu\in\hfrak^*$ with $\lambda\preceq \nu$.    Then there exist a $\gfrak$-module filtration $0=M_0\subseteq M_1\subseteq M_2\subseteq \ldots \subseteq M_{k-1}\subseteq M_k =M$ and a subset $J\subseteq \{1,2,\ldots,k\}$ such that
	\begin{enumerate}[(i)]
		\item if $j\in J$, then $M_j/M_{j-1}\cong \Llie\big(\xi(j)\big)$ for some $\xi(j)\in\hfrak^*$ with $\lambda\preceq \xi(j)\preceq \nu$,
		\item if $j\notin J$, then either $\left(M_j/M_{j-1}\right)^\mu=0$ holds for every $\mu\in\hfrak^*$ satisfying $\lambda\preceq\mu\preceq\nu$, or $M_j/M_{j-1}\cong \Llie\big(\xi(j)\big)$ for some $\xi(j)\succ\lambda$ such that $\xi(j)\not\preceq \nu$.
	\end{enumerate}	Such a filtration is called a \emph{composition series of $M$ with bounds $\lambda$ and $\nu$}.  The set $J$ is called the \emph{relevant index set} of such a filtration.
	\label{cor:kac_series}
	
	\begin{pf}
		Since the interval $[\lambda,\nu]:=\big\{\zeta\in\hfrak^*\,\suchthat{\lambda\preceq\zeta\preceq \nu}\big\}$ is finite (as $\bfrak$ is Dynkin) and $M$ is $\nfrak$-locally finite, the submodule
		\begin{align}
			\tilde{M}:=\sum_{\substack{{\mu\in\hfrak^*}\\{\lambda\preceq \mu\preceq\nu}}}\,\Ulie(\gfrak)\cdot M^\mu =\sum_{\substack{{\mu\in\hfrak^*}\\{\lambda\preceq \mu\preceq\nu}}}\,\Ulie(\nfrak)\cdot M^\mu
			\label{eq:mtilde}
		\end{align}
		has finitely many weights $\zeta$ with $\zeta\succeq \lambda$.  Therefore, $\tilde{M}$ has an upper bound $\upsilon\in\hfrak^*$.  We apply Lemma~\ref{lem:kac_series} on $\tilde{M}$ and obtain a filtration
		\begin{align}
			0=M_0\subseteq M_1\subseteq M_2\subseteq \ldots\subseteq M_{k-2}\subseteq M_{k-1}=\tilde{M}
		\end{align}
		along with a subset $\tilde{J}\subseteq\{1,2,\ldots,k-1\}$ satisfying the conditions that, for every element $j\in\tilde{J}$, $M_j/M_{j-1}\cong \Llie\big(\xi(j)\big)$ for some $\xi(j)\in\hfrak^*$ with $\xi(j)\succeq \lambda$, and that, whenever $j\notin\tilde{J}$, $\left(M_j/M_{j-1}\right)^\mu=0$ for every $\mu\succeq\lambda$.  Then, by setting $M_{k}:=M$, we have the filtration 
		\begin{align}
			0=M_0\subseteq M_1\subseteq M_2\subseteq \ldots\subseteq M_{k-2}\subseteq M_{k-1}\subseteq M_{k}=M\,.
			\label{eq:filtrationM}
		\end{align}
		Let $J:=\big\{j\in\tilde{J}\,|\,\xi(j)\preceq \nu\big\}$.  The filtration (\ref{eq:filtrationM}) clearly satisfies (i) and (ii), with the relevant index set $J$, noting that $\left(M_k/M_{k-1}\right)^\mu=\left(M/\tilde{M}\right)^\mu=0$ for all $\mu\in\hfrak^*$ such that $\lambda\preceq\mu\preceq \nu$ holds.
	\end{pf}
\end{cor}

\begin{define}
	Let $\lambda,\nu\in\hfrak^*$ satisfy $\lambda\preceq \nu$.  Suppose that 
	\begin{align} 0=M_0\subseteq M_1\subseteq M_2\subseteq \ldots \subseteq M_{k-1}\subseteq M_k=M\label{eq:com1}\end{align} and \begin{align} 0=M'_0\subseteq M'_1\subseteq M'_2\subseteq \ldots \subseteq M'_{k'-1}\subseteq M'_{k'}=M\label{eq:com2}\end{align} are two composition series of $M\in\bar\bggO$ with bounds $\lambda$ and $\nu$, and with relevant index sets $J$ and $J'$, respectively.  We say that these filtrations are \emph{equivalent} if there exists a bijection $f:J\to J'$ such that
	$
		M_{j}/M_{j-1}\cong M'_{f(j)}/M'_{f(j)-1}$ for all $j\in J$.
\end{define}

\begin{lem}
	Let $\lambda,\nu\in\hfrak^*$ be such that $\lambda\preceq \nu$.  Denote by $\tilde{M}$ the submodule of $M$ given by (\ref{eq:mtilde}).  Suppose that $0=M_0\subseteq M_1 \subseteq M_2 \subseteq \ldots \subseteq M_{k-1}\subseteq M_k=M$.  Define $\tilde{M}_j:=M_j\cap\tilde{M}$ for every $j=0,1,2,\ldots,k$.  Then, $0=\tilde{M}_0\subseteq \tilde{M}_1\subseteq \tilde{M}_2\subseteq \ldots\subseteq \tilde{M}_{k-1}\subseteq \tilde{M}_k=\tilde{M}$ is a composition series of $\tilde{M}$ with bounds $\lambda$ and $\nu$.
	
	\begin{pf}
		For each $j=1,2,\ldots,k$, define $\varphi_j:\tilde{M}_j/\tilde{M}_{j-1}\to M_j/M_{j-1}$ via $v+\tilde{M}_{j-1}\mapsto v+M_{j-1}$ for every $v\in\tilde{M}_j$.  Clearly, $\varphi_j$ is well defined, and it is injective because 
		\begin{align} 	\tilde{M}_{j}\cap M_{j-1}=\tilde{M}\cap M_j\cap M_{j-1}=\tilde{M}\cap M_{j-1}=\tilde{M}_{j-1}\,.\end{align}
		If $M_j/M_{j-1}\cong \Llie(\mu)$ for some $\mu\in\hfrak^*$ with $\mu\preceq\lambda$, then 
		\begin{align} \dim_\amsbb{K}\left(\tilde{M}_j^\mu\right)=\dim_\amsbb{K}\left(M_j^\mu\right)=\dim_\amsbb{K}\left(M_{j-1}^\mu\right)+1=\dim_\amsbb{K}\left(\tilde{M}_{j-1}^\mu\right)+1\,.\end{align}
		Hence, 
		$
			\dim_\amsbb{K}\left(\left(\tilde{M}_{j}/\tilde{M}_{j-1}\right)^\mu\right)=1$, so $\varphi_j$ is nonzero.  As $M_j/M_{j-1}$ is simple, $\varphi_j$ must be surjective, whence it gives an isomorphism 
		$\tilde{M}_j/\tilde{M}_{j-1} \cong M_j/M_{j-1}\cong \Llie(\mu)$.
		
		Let $J$ be the relevant index set of the composition series 
		$
			M_0\subseteq M_1 \subseteq \ldots \subseteq M_{k-1}\subseteq M_k
		$
		of $M$ with bounds $\lambda$ and $\nu$.  By the observation above, if $j\in J$, then $\tilde{M}_j/\tilde{M}_{j-1}\cong M_j/M_{j-1}$ is simple with highest weight $\mu$ satisfying $\lambda\preceq\mu\preceq\nu$. If an index $j\in\{1,2,\ldots,k\}\smallsetminus J$ is such that $M_j/M_{j-1}$ is a simple module with highest weight $\xi\succ \lambda$ with $\xi\not\preceq\lambda$, then using the embedding  $\varphi_j: \tilde{M}_j/\tilde{M}_{j-1}\to M_j/M_{j-1}$, we  conclude that either 
		$
			\tilde{M}_j/\tilde{M}_{j-1}\cong\Llie(\xi)$ or $\tilde{M}_j/\tilde{M}_{j-1}=0$.  Finally, if $j\in\{1,2,\ldots,k\}\smallsetminus J$ is such that $\left(M_j/M_{j-1}\right)^\mu =0$ for every $\mu$ with $\lambda\preceq\mu\preceq\nu$, using the embedding $\varphi_j: \tilde{M}_j/\tilde{M}_{j-1}\to M_j/M_{j-1}$, we see that 
			$
				\left(\tilde{M}_j/\tilde{M}_{j-1}\right)^\mu =0$
		 for every $\mu\in\hfrak^*$ with $\lambda\preceq\mu\preceq\nu$.
	\end{pf}
\end{lem}

\begin{thm}
	Let $\lambda,\nu\in\hfrak^*$ be such that $\lambda\preceq \nu$.  Then, any two composition series of $M\in\bar\bggO$ with bounds $\lambda$ and $\nu$ are equivalent.
	\label{thm:kacequiv}
	
	\begin{pf}
		Suppose that (\ref{eq:com1}) and (\ref{eq:com2}) are two composition series of $M$ with bounds $\lambda$ and $\nu$.  Let $\tilde{M}$ be the submodule of $M$ defined by (\ref{eq:mtilde}).  From the lemma above, it suffices to assume that $M=\tilde{M}$.  
		
		From the assumption $M=\tilde{M}$, there are finitely many weights $\mu$ of $M$ satisfying $\mu\preceq\lambda$.  Thus, we can refine (\ref{eq:com1}) and (\ref{eq:com2}) in the same manner as in Theorem~\ref{lem:kac_series} to get index sets $\tilde{J}\subseteq\{1,2,\ldots,k\}$ and $\tilde{J}'\subseteq\{1,2,\ldots,k'\}$ such that the following three conditions are met:
		\begin{enumerate}[(i)]
			\item $J\subseteq \tilde{J}$ and $J'\subseteq \tilde{J}'$,
			\item for $j\in \tilde{J}$ and $j'\in \tilde{J}'$, $M_j/M_{j-1}$ and $M'_{j'}/M'_{j'-1}$ are simple modules with highest weights greater than or equal to $\lambda$,
			\item for $j\notin\tilde{J}$ and $j'\notin\tilde{J}'$, all the weight spaces $\left(M_j/M_{j-1}\right)^\mu$ and $\left(M'_{j'}/M'_{j'-1}\right)^\mu$ with $\mu\succeq\lambda$ are the zero vector space.
		\end{enumerate}
		As a result, we can instead show that there exists a bijective function $\tilde{f}:\tilde{J}\to\tilde{J}'$ such that $M_j/M_{j-1}\cong M'_{f(j)}/M'_{f(j)-1}$ for all $j\in\tilde{J}$.  The restriction $f:=\tilde{f}|_J$ yields a bijection $f:J\to J'$ as required.
		
		For each $\xi\in\hfrak^*$ with $\lambda \preceq \xi \preceq \nu$, let $t(\xi)$ denote the maximum possible value of the length of the positive root $\mu-\xi$, where $\mu\succeq\xi$ is a singular weight of $M$.  For each $l=0,1,2,\ldots$, write $T_l$ for the set $\big\{\xi\in\hfrak^*\suchthat{}\xi\succeq\lambda\textup{ and }t(\xi)=l\big\}$.  We shall instead prove that, for a fixed $l=0,1,2,\ldots$, the number of $j\in J$ with $M_{j}/M_{j-1}\cong \Llie(\xi)$ is the same as the number of $j'\in J'$ with $M'_{j'}/M'_{j'-1}\cong \Llie(\xi)$ for every $\xi \in T_l$.
		
		The proof goes by induction on $l$.  For the base case $l=0$, every $\xi\in T_l$ is a singular weight of $M$, whence the weight space $M^\xi$ comes from  $\dim_\amsbb{K}\left(M^\xi\right)$ copies of $\Llie(\xi)$ in any composition series with  bounds $\lambda$ and $\nu$.   
		
		Let now assume that $l>0$ and $\xi\in T_l$.  By the induction hypothesis, the multiplicities of each factor $\Llie(\tilde{\xi})$ with $\tilde{\xi}\in T_0\cup T_1 \cup\ldots \cup T_{l-1}$ in the filtrations (\ref{eq:com1}) and (\ref{eq:com2}) are equal, and let $m\left(\tilde{\xi}\right)$ denote the common value.  For each $j\in J$ with $M_j/M_{j-1} \not\cong \Llie\left(\tilde{\xi}\right)$ with $\tilde{\xi}\in T_0\cup T_1\cup\ldots\cup T_{l-1}$, we observe that either $\left(M_j/M_{j-1}\right)^{\xi}=0$, or $\xi$ is a singular weight of $M_j/M_{j-1}$ (making $M_j/M_{j-1}\cong \Llie(\xi)$).  Hence, there are exactly 
		\begin{align}
		m(\xi):= \dim_\amsbb{K}\left(M^\xi\right)-\sum_{r=0}^{l-1}\,\sum_{\tilde{\xi}\in T_r}\,m(\tilde{\xi})\,\dim_\amsbb{K}\Big(\big(\Llie(\tilde{\xi})\big)^\xi\Big)
		\label{eq:multcount}
		\end{align} values of such $j\in J$ with $M_j/M_{j-1}\cong \Llie(\xi)$.  Therefore, $m(\xi)$ is the multiplicity of $\Llie(\xi)$ in (\ref{eq:com1}).  Since the value $m(\xi)$ as shown in (\ref{eq:multcount}) depends only on the previously known values $m\left(\tilde{\xi}\right)$ with $\tilde{\xi}\in T_0\cup T_1 \cup\ldots \cup T_{l-1}$, $m(\xi)$ is also the multiplicity of $\Llie(\xi)$ in (\ref{eq:com2}).  The induction is now complete and the claim follows.
	\end{pf}
\end{thm}

Corollary~\ref{cor:kac_series} and Theorem~\ref{thm:kacequiv} form a partial extension of the usual Jordan-H\"{o}lder Theorem for modules of finite length.  Based on this, we now extend the usual definition of composition factors and composition factor multiplicities as follows.

\begin{cor}
	Let $M\in\bar\bggO$ and $\mu\in\hfrak^*$ be such that $\lambda\preceq\mu\preceq \nu$.  The number of times $\Llie(\mu)$ occurs as a factor in any composition series of $M$ with bounds $\lambda$ and $\nu$ is independent of the choice of the composition series with bounds and the choice of the bounds $\lambda,\nu\in\hfrak^*$ (as long as $\lambda\preceq \mu \preceq \nu$).    This number is known as the \emph{composition factor multiplicity} of $\Llie(\mu)$ in $M$, and is denoted by $\big[M:\Llie(\mu)\big]$.  If $\big[M:\Llie(\mu)\big]>0$, then we say that $\Llie(\mu)$ is a \emph{composition factor} of $M$.
	\label{cor:compositionfactormult}
	
	\begin{pf}
		For given weights $\lambda,\nu\in\hfrak^*$, Theorem~\ref{thm:kacequiv} guarantees that the number $m^M_\mu(\lambda,\nu)$  of times $\Llie(\mu)$ occurs as a factor does not depend on the choice of the composition series of $M$ with bounds $\lambda$ and $\nu$.  We have to show that $m^M_\mu(\lambda,\nu)$ is also independent of $\lambda$ and $\nu$, provided that $\lambda\preceq\mu\preceq\nu$.
		
		Let $\lambda,\nu\in\hfrak^*$ be such that $\lambda\preceq\mu\preceq\nu$.  We choose an arbitrary composition series $0=M_0\subsetneq M_1\subsetneq M_2\subsetneq \ldots \subsetneq M_k=M$ of $M$ with bounds $\lambda$ and $\nu$.  Then, this filtration is also a composition series with bounds $\mu$ and $\mu$.  Again, by Theorem~\ref{thm:kacequiv}, this filtration is equivalent to any composition series with bounds $\mu$ and $\mu$, which immediately implies that 
		$m^M_\mu(\lambda,\nu)=m^M_\mu(\mu,\mu)$.
	\end{pf}
\end{cor}

Now we shall use the composition series with bounds to study generalized composition series, as introduced below.  With the restriction that the modules in $\bbggO$ have finite-dimensional weight spaces, we shall see that these generalized composition series behave similarly to the composition series of modules of finite length. 

\begin{define}
	A \emph{generalized composition series} of a module $M\in\bbggO$ is a family of submodules $\left(M_j\right)_{j\in \mathcal{J}}$ satisfying the following conditions:
	\begin{enumerate}[(i)]
		\item the index set $\mathcal{J}$ is equipped with a total order $<$,
		\item $M_j\subsetneq M_k$ for each $j,k\in \mathcal{J}$ with $j<k$.
		\item $\bigcap_{j\in\mathcal{J}}\,M_{j}=0$ and $\bigcup_{j\in\mathcal{J}}\,M_j=M$,
		\item $M_j/M_{<j}$ is a simple module for all $j\in \mathcal{J}$, where $M_{<j}:=\bigcup_{k<j}\,M_k$.
	\end{enumerate}
	\label{def:gencompseries}
\end{define}
%

\begin{thm}
	Every $M\in\bar\bggO$ has a generalized composition series. 
	\label{thm:gencomp} 
	
	\begin{pf}
		First, we shall prove this theorem when $M$ is indecomposable.  We start with arbitrary weights $\lambda(1)$ and $\nu(1)$ of $M$ with $\lambda(1)\preceq \nu(1)$.  Let $\mathcal{J}(1)$ be the relevant index set of a composition series of $M$ with bounds $\lambda(1)$ and $\nu(1)$.  We create two sequences of weights $\big\{\lambda(k)\big\}_{k\in\mathbb{Z}_{>0}}$ and $\big\{\nu(k)\big\}_{k\in\mathbb{Z}_{>0}}$ in such a way that
		\begin{align}
			\lambda(1)\succ\lambda(2)\succ\lambda(3)\succ\ldots\,,
			\label{eq:lambdakk}
		\end{align}
		\begin{align}
			\nu(1)\prec\nu(2)\prec\nu(3)\prec\ldots\,,
			\label{eq:nukk}
		\end{align}
		and, for every weight $\zeta\in \textup{supp}(M)$, there exists $k\in\amsbb{Z}_{>0}$ (depending on $\zeta$) such that $\lambda(k)\preceq \zeta\preceq \nu(k)$.  Note that $\lambda(k)$ and $\nu(k)$ do not have to be weights of $M$.  Therefore, it is always possible to find an infinite strictly decreasing sequence (\ref{eq:lambdakk}) and an infinite strictly increasing sequence (\ref{eq:nukk}).

		Suppose that $\mathcal{J}(k)$ is known.  We extend the filtration in the $k$-th step to obtain a composition series of $M$ with bounds $\lambda(k+1)$ and $\mu(k+1)$.  To be precise, suppose that
		\begin{align}
			0=M^k_0\subsetneq M^k_1\subsetneq M^k_2 \subsetneq \ldots \subsetneq M^k_{l(k)-1}\subsetneq M^k_{l(k)}=M
			\label{eq:lambdak}
		\end{align}
		is a composition series with bounds $\lambda(k)$ and $\mu(k)$.  If $i$ is in the relevant index set $\mathcal{J}(k)$, then $M^k_{i}/M^k_{i-1}$ is simple with highest weight $\mu$ with 
		$
			\lambda(k+1)\prec\lambda(k)\preceq \mu\preceq \nu(k)\prec \nu(k+1)$.
		If $i>0$ is not in the relevant index set, then we find a composition series of $M^k_i/M^k_{i-1}$ with bounds $\lambda(k+1)$ and $\nu(k+1)$.  Then, take the preimages of the submodules that occur this composition series of $M^k_i/M^k_{i-1}$ for each $i>0$ not in $\mathcal{J}(k)$.  Using these preimages, we then refine the composition series (\ref{eq:lambdak}) and obtain a composition series 
		\begin{align}
			0=M^{k+1}_0\subsetneq M^{k+1}_1\subsetneq M^{k+1}_2 \subsetneq \ldots \subsetneq M^{k+1}_{l(k+1)-1}\subsetneq M^{k+1}_{l(k+1)}=M
		\end{align}
		with bounds $\lambda(k+1)$ and $\nu(k+1)$, along with an inclusion $\iota_k:\mathcal{J}(k)\to\mathcal{J}(k+1)$ of totally ordered sets satisfying $M^k_i/M^k_{i-1}\cong M^{k+1}_{\iota_k(i)}/M^{k+1}_{\iota_k(i)-1}$ for every $j\in\mathcal{J}(k)$.

		We now take $\mathcal{J}$ to be the direct limit $\lim\limits_{\underset{k}{\boldsymbol{\longrightarrow}}}\,\mathcal{J}(k)$.  By construction, there is a total order $<$ on $\mathcal{J}$ extending the total orders on the sets $\mathcal{J}(k)$.  Each $j\in\mathcal{J}$ corresponds to an element in $j(k)\in \mathcal{J}(k)$ for some large enough $k$, and to a submodule $M_j:=M^k_{j(k)}$ of $M$ in the composition series from the $k$-th step.  Note that 
		$M_{<j}=\bigcup_{j'<j}\,M_j=M^k_{j(k)-1}$, whence $M_j/M_{<j}=M^k_{j(k)}/M^k_{j(k)-1}$ is simple.   Clearly, the index set $\mathcal{J}$ and  the family of submodules $\left(M_j\right)_{j\in\mathcal{J}}$ satisfy the requirements.
		
		Finally, suppose that $M$ has a direct sum decomposition $M=\bigoplus_{t\in \pazocal{I}}\,D_t$, where each $D_t$ is indecomposable (by Theorem~\ref{thm:directsumdecomposition}).  We first equip $\pazocal{I}$ with a well order $\triangleleft$ (which exists by the Well-Ordering Principle).  Then, we create a generalized filtration series $\left\{D_t(j)\right\}_{j\in \mathcal{J}_t}$ for each $D_t$.  Write $<_t$ for the total order on $\mathcal{J}_t$.  Let $\pazocal{J}$ be the totally ordered set
		$
		\big\{(t,j)\,\suchthat{t\in\pazocal{I}\textup{ and }j\in \mathcal{J}_t}\big\}$
		with the total order $<$ defined via the lexicographic ordering as follows:
		\begin{align}
			(t,j)<(\tilde{t},\tilde{j})\textup{ if and only if }t\triangleleft \tilde{t}\,,\textup{ or }t=\tilde{t}\textup{ and }j\,{<_t}\,\tilde{j}\,.
		\end{align}
		Then, we take $M_{(t,j)}:=D_t(j)\oplus \left(\bigoplus_{\tilde{t}\triangleleft t}\,D_{\tilde{t}}\right)$ for every $(t,j)\in\pazocal{I}$.  Then, it is obvious that $\left\{M_{(t,j)}\right\}_{(t,j)\in \pazocal{I}}$ is a generalized composition series of $M$.
	\end{pf}
	
\end{thm}

\begin{define}
	\label{def:equivcomp}
	Two generalized composition series $\left(M_j\right)_{j\in\mathcal{J}}$ and $\left(M'_{j'}\right)_{j'\in\mathcal{J}'}$ are \emph{equivalent} if there exist a bijective function $f:\mathcal{J}\to\mathcal{J}'$ such that 
	$
		M_j/M_{<\,j}\cong M'_{f(j)}/M'_{<'\,f(j)}
	$ for each $j\in\mathcal{J}$.  Here, $<$ is the total order on $\mathcal{J}$, whereas $<'$ is the total order on $\mathcal{J}'$.  In addition, $M_{< \,j}:=\bigcup_{k\,<\, j}\,M_k$ as well as $M'_{<'\,j'}:=\bigcup_{k'\,<'\,j'}\,M'_{k'}$ for all $j\in\mathcal{J}$ and $j'\in\mathcal{J}'$.
\end{define}


\begin{thm}
	Any two generalized composition series of $M\in\bar\bggO$ are equivalent.  
	\label{thm:uniquegencomp}
	
	\begin{pf}
		We may first assume that $M$ is indecomposable.  Let $\left(M_j\right)_{j\in \mathcal{J}}$ and $\left(M'_{j'}\right)_{j'\in\mathcal{J}'}$ be two generalized composition series of an object $M\in\bbggO$. We create a decreasing sequence of weights $\big\{\lambda(k)\big\}_{k\in\mathbb{Z}_{>0}}$ and an increasing sequence of weights $\big\{\nu(k)\big\}_{k\in\mathbb{Z}_{>0}}$ such that every weight $\mu \in\textup{supp}(M)$ satisfies $\lambda(k)\preceq\mu\preceq\nu(k)$ for some $k$. 
		
		For each $k$, define
		\begin{align}
			\mathcal{J}(k):= \Big\{j\in \mathcal{J}\,\suchthat{M_j/M_{<\,j}\cong \Llie(\xi)\textup{ for some }\xi\textup{ with }\lambda(k)\preceq\xi\preceq\nu(k)}\Big\}
		\end{align}
		and
		\begin{align}			
			\mathcal{J}'(k):= \Big\{j'\in \mathcal{J}'\,\suchthat{M'_j/M'_{<'\,j'}\cong \Llie(\xi)\textup{ for some }\xi\textup{ with }\lambda(k)\preceq\xi\preceq\nu(k)}\Big\}
		\end{align}
		From Theorem~\ref{thm:kacequiv}, we have a bijection $f_k:\mathcal{J}(k)\to\mathcal{J}'(k)$ satisfying the condition that $M_j/M_{<\,j}\cong M'_{f_k(j)}/M'_{<'\,f_k(j)}$ for every $j\in\mathcal{J}(k)$.  We claim that the bijections $f_k$ can be aligned so that $\left.f_{k+1}\right|_{\mathcal{J}(k)}=f_k$.   Taking the direct limit $f:=\lim\limits_{\underset{k}{\boldsymbol{\longrightarrow}}}\,f_k$ then yields a bijection $f:\mathcal{J}\to\mathcal{J}'$ satisfying the requirement of Definition~\ref{def:equivcomp}.
		
		To prove the claim above, assume that $\left.f_{k+1}\right|_{\mathcal{J}(k)}\neq f(k)$.  Then we define the function $\tilde{f}_{k+1}:\mathcal{J}(k+1)\to\mathcal{J}'(k+1)$ as follows: 
		\begin{align}
			\tilde{f}_{k+1}(j)=\left\{\begin{array}{ll}
				f_k(j) &\textup{if }j\in \mathcal{J}(k)\,,\\
				\tilde{f}_{k+1}(j)&	\textup{if }j\in\mathcal{J}(k+1)\smallsetminus\mathcal{J}(k)\,.
			\end{array}\right.
		\end{align}  
		Replacing $f_{k+1}$ by $\tilde{f}_{k+1}$ and continuing this process for each positive integer $k$, we obtain a set of aligned bijections as desired.

		When $M$ is decomposable, we note that it has a unique direct sum decomposition into indecomposable direct summands (Theorem~\ref{thm:directsumdecomposition} and Corollary~\ref{cor:KSRA}).  From this, we can easily conclude that two generalized composition series of $M$ are equivalent.
	\end{pf}
\end{thm}


\begin{define}
\label{def:Wlambda}
Let $\lambda\in\hfrak^*$.  Let $W[\lambda]$ be the subgroup of $W$ containing all $w\in W$ such that $w\cdot\lambda-\lambda \in \Lambda$.  Write $W_n$ for the Weyl group of $\gfrak_n$.    We similarly define $W_n\left[\lambda_n\right]$ for each $n\in\amsbb{Z}_{>0}$ and $\lambda_n\in\hfrak_n^*$.  These subgroups are known as the \emph{integral Weyl groups}.
\end{define}

\begin{thm}
	Let $\lambda,\mu\in\hfrak^*$.  If the simple module $\Llie(\mu)$ is a composition factor of the Verma module $\Mlie(\lambda)$, then $\mu\preceq \lambda$ and $\mu \in W[\lambda]\cdot\lambda$.
	\label{thm:compositionfactors}
	
	\begin{pf}
		Suppose that $\Mlie(\lambda)$ has $\Llie(\mu)$ as a composition factor.   Let $v$ be a highest-weight vector of $\Mlie(\lambda)$.  Then, for all sufficiently large $n\in\mathbb{N}$, $M_n:=\Ulie\left(\gfrak_n\right)\cdot v \in\bggO^{\gfrak_n}_{\hfrak_n}$ is a Verma module and must have $\gfrak_n$-submodules $N_n$ and $N'_n$ with $N_n\subseteq N'_n$ and $N'_n/N_n$ has $\mu_n:=\mu|_{\hfrak_n}$ as a highest weight.  Hence, $\Llie\left(\mu_n\right)$ is a composition factor of $M_n\cong\Mlie\left(\lambda_n\right)$, where $\lambda_n:=\lambda|_{\hfrak_n}$.  Due to the finite-dimensional theory, 
		$
			\mu_n\preceq\lambda_n\textup{ and }\mu_n\in W_n\left[\lambda_n\right]\cdot\lambda_n$.  The result follows immediately.
	\end{pf}
\end{thm}

\section{Block Decomposition and Kazhdan-Lusztig Multiplicities}
\label{ch:blocks}

In this section, $\gfrak$ is a root-reductive Lie algebra and $\bfrak$ is a Dynkin Borel subalgebra.  Our objective is to prove that the block decomposition of $\bbggO$ is similar to that of the category $\bggO$. 

\subsection{Block Decomposition}

Let $M$ be an indecomposable object of the category $\bbggO$.  We shall first construct a countable ordered set $\Gamma(M)=\left(v_0,v_1,v_2,\ldots\right)$ such that $\Gamma(M)\subseteq M$ consists of weight vectors of $M$ which generate $M$ as a $\Ulie(\gfrak)$-module, and the set $\Gamma(M)$ has certain desirable properties.  If $M=0$, then we set $\Gamma(M):=(0,0,0,0,0,\ldots)$.

For $M\neq 0$, we let $u\neq 0$ be a singular vector of $M$, and $\xi\in\hfrak^*$ the weight of $u$.  Let $[\xi]$ denote the set of all weights $\zeta\in\hfrak^*$ such that $\zeta-\xi$ is in the $\amsbb{Z}$-span of the simple ($\bfrak$-positive) roots of $\gfrak$.  For a weight $\zeta\in[\xi]$, the \emph{distance} between $\zeta$ and $\xi$, denoted by $\textup{dist}(\zeta,\xi)$, is defined to be the sum 
\begin{align}
	\sum_{i=1}^k\,\left|t_i\right|\,,\textup{ if }\zeta-\xi=\sum_{i=1}^k\,t_i\,\alpha_i\,,
\end{align}
where $\alpha_1,\alpha_2,\ldots,\alpha_k$ are ($\bfrak$-positive) simple roots of $\gfrak$ and $t_1,t_2,\ldots,t_k\in\amsbb{Z}$.  Furthermore, the height of $\zeta-\xi$, denoted by $\textup{ht}(\zeta-\xi)$, is the smallest integer $n\geq 0$ such that $h_{\alpha_1},h_{\alpha_2},\ldots,h_{\alpha_k}$ are all in $\gfrak_n$ (noting that $\textup{ht}(\zeta-\xi)=0$ if and only if $\zeta =\xi$).  

We start with $m:=0$; then we set $d(0):=0$ and $v_0:=u$. Now, for $m>0$, suppose that the value $d(m-1)$ is known and that the vectors $v_0,v_1,\ldots,v_{d(m-1)}$ have been defined.  The set $\pazocal{S}_m$ of weights $\zeta \in [\xi]$ such that $\textup{dist}(\zeta,\xi)\leq m$ and $\textup{ht}(\zeta-\xi)\leq m$ is a finite set.  Let $\pazocal{V}^1_m$  denote the $\amsbb{K}$-span of all weight vectors $v\in M$ with weights in $\pazocal{S}_m$ such that $\nfrak\cdot v=0$.  Let $u^1_1,u^1_2,\ldots,u^1_{l_1}$ be weight vectors of $M$ which form a $\amsbb{K}$-basis of $\pazocal{V}^1_m$.  

Assume that the collections $\left(u^1_j\right)_{j=1}^{l_1}$, $\left(u^2_j\right)_{j=1}^{l_2}$, $\ldots$, $\left(u^{r}_j\right)_{j=1}^{l_{r}}$ of weight vectors of $M$ have been obtained. Consider the module 
\begin{align}
	\underline{M}(m,r):=M/\left(\sum_{i=0}^{d(m-1)}\,\Ulie(\gfrak)\cdot v_i+ \sum_{p=1}^r\,\sum_{j=1}^{l_p}\,\Ulie(\gfrak)\cdot u^p_j\right)\,.
\end{align}
Let $\pazocal{V}^{r+1}_m$ denote the $\amsbb{K}$-span of all weight vectors $v\in \underline{M}(m,r)$ with weights in $\pazocal{S}_m$ such that $\nfrak\cdot v=0$.  Take $\tilde{u}^{r+1}_1,\tilde{u}^{r+1}_2,\ldots,\tilde{u}^{r+1}_{l_{r+1}}$ to be weight vectors of $M$ which form a $\amsbb{K}$-basis of $\pazocal{V}^{r+1}_m$.  Now, there exist weight vectors $u^{r+1}_1,u^{r+1}_2,\ldots,u^{r+1}_{l_{r+1}}$, whose respective images under the projection $M\to \underline{M}(m,r)$ are $\tilde{u}^{r+1}_1,\tilde{u}^{r+1}_2,\ldots,\tilde{u}^{r+1}_{l_{r+1}}$.

The process in the previous paragraph must end with $\pazocal{V}_m^{\bar{r}(m)+1}=0$ for some nonnegative integer $\bar{r}(m)$ because the vector subspace of $M$ spanned by the weight vectors with weights in $\pazocal{S}_m$ is finite dimensional.  Then, we take 
\begin{align}
	d(m):=d(m-1)+\sum_{p=1}^{\bar{r}(m)}\,l_p\,,
\textup{ and }
	v_{d(m-1)+j}:=u^{l_p}_{j-\sum_{\tilde{p}=1}^{p-1}\,l_{\tilde{p}}}
\end{align}
if $\sum_{\tilde{p}=1}^{p-1}\,l_{\tilde{p}}<j\leq \sum_{\tilde{p}=1}^p\,l_{\tilde{p}}$.

Note that $d(m)>d(m-1)$ for every $m=1,2,\ldots$ because $\pazocal{V}_m^1$ always contains $u$.  When $M$ is a $\gfrak$-module of finite length, it is possible that $\Gamma(M)$ is eventually periodic (that is, there exist positive integers $n_0$ and $n_1$ such that $v_{n}=v_{n+n_1}$ for every integer $n\geq n_0$).  In particular, the ordered set $\Gamma(M)$ may take the form $(u,u,u,u,\ldots)$ when $M$ is a highest-weight module with $u$ as a highest-weight vector. 

 We claim that the ordered set 
 $\Gamma(M):=\left(v_0,v_1,v_2,v_3,\ldots\right)$ generates $M$ as a $\Ulie(\gfrak)$-module.  For a fixed weight $\zeta$ of $M$, consider the vector subspace $M^\zeta$.  

Let $\pazocal{T}_\zeta$ denote the set of all weights $\tilde{\zeta}$ of the $\Ulie(\gfrak)$-module $\Ulie(\gfrak)\cdot M^\zeta$ which satisfy $\tilde{\zeta}\succeq \zeta$.  Note that $\pazocal{T}_\zeta$ is finite as $\nfrak$ acts locally finitely on $M^\zeta$.  Let $m_\zeta$ denote the maximum value of the two numbers $
		\max\big\{\textup{dist}(\tilde{\zeta},\xi)\suchthat{}\tilde{\zeta} \in \pazocal{T}_\zeta\big\}$ and $\max\big\{\textup{ht}(\tilde{\zeta}-\xi)\suchthat{}\tilde{\zeta} \in \pazocal{T}_\zeta\big\}$.
Then, in the $m_\zeta$-th step of the procedure (from which $d\left(m_\zeta\right)$ is obtained), the $\Ulie(\gfrak)$-module
\begin{align}
	\underline{M}\Big(m_\zeta,\bar{r}\left(m_\zeta\right)\Big)=M/\left(\sum_{i=0}^{d(m-1)}\,\Ulie(\gfrak)\cdot v_i + \sum_{p=1}^{\bar{r}(m)}\,\sum_{j=1}^{l_p}\,\Ulie(\gfrak)\cdot u^p_j\right)
\end{align}
cannot have $\Llie(\zeta)$ as a composition factor.  To elaborate, if such a composition factor exists, it must arise from the weight space of $\underline{M}\Big(m_\zeta,\bar{r}\left(m_\zeta\right)\Big)$ with weight $\zeta$.  However, all composition factors of $M$ isomorphic to $\Llie(\zeta)$ come from subquotients of  $\Ulie(\gfrak)\cdot M^\zeta$, and by the definition of $m_\zeta$, the image of $\Ulie(\gfrak)\cdot M^\zeta$ under the canonical projection $M\to \underline{M}\Big(m_\zeta,\bar{r}\left(m_\zeta\right)\Big)$ has no composition factors isomorphic to $\Llie(\zeta)$.

Therefore, for every composition factor $\Llie(\zeta)$ of $M$, it is exhausted in the quotient module $\underline{M}\Big(m_\zeta,\bar{r}\left(m_\zeta\right)\Big)$.  Thus, the sum $\sum_{i=0}^\infty\,\Ulie\left(\gfrak\right)\cdot v_i$ is precisely $M$ itself.  

It is important to note the following properties of the set $\Gamma(M)$.  For every $n\in\amsbb{Z}_{>0}$, let $\Xi_n\subseteq \hfrak^*$ be the support (as a semisimple $\hfrak$-module) of the $\gfrak'_n$-module $M_n:=\sum_{i=0}^n\,\Ulie\left(\gfrak'_n\right)\cdot v_i$, where $\gfrak_n':=\hfrak+\gfrak_n$.   Then, for every $\xi\in \Xi_n$ and for any integer $\tilde{n}\geq n$, we have
	\begin{align}
		\dim_\amsbb{K}\left(M_n^{\xi}\right)=\dim_\amsbb{K}\Big(\big(\Ulie\left(\gfrak'_{\tilde{n}}\right)\cdot M_n\big)^{\xi}\Big)\,.
		\label{eq:equaldim}
	\end{align}
	This is because our construction of $\Gamma(M)$ ensures that 
	$\Ulie\left(\gfrak'_{\tilde{n}}\right)\cdot M_n = \Ulie\left(\nfrak_{\tilde{n}}^-\right)\cdot M_n$.  Note that, if $x$ is in a positive root space of $\gfrak'_{\tilde{n}}$ that is not in $\gfrak'_n$, then 
	$x\cdot M_n=0$.

	\begin{define}
		Let $\lambda\in\hfrak^*$.  Define $\bbggO\llbracket\lambda\rrbracket$ to be the full subcategory of $\bbggO$ consisting of modules $M$ whose composition factors are of the form $\Llie(\mu)$ with $\mu\in\llbracket\lambda\rrbracket$, where $\llbracket\lambda\rrbracket$ is the integral Weyl dot-orbit $\llbracket\lambda\rrbracket:=W[\lambda]\cdot\lambda$.
	\end{define}
	
	\begin{prop}
		Let $M\in\bbggO$ be indecomposable and $\lambda\in\hfrak^*$ be such that $\Llie(\lambda)$ is a composition factor of $M$.  Then, all composition factors of $M$ are of the form $\Llie(\mu)$ for some $\mu\in \llbracket\lambda\rrbracket$.
		\label{prop:blocks}
		
		\begin{pf}
			Let $\lambda$ and $\mu$ be on different integral Weyl dot-orbits.  Suppose there exists an indecomposable $M\in\bar\bggO$ with $\Llie(\lambda)$ and $\Llie(\mu)$ as composition factors.  Since $M\neq0$, we can apply the algorithm discussed earlier in this subsection and obtain an ordered set $\Gamma(M)=\left(v_0,v_1,v_2,\ldots\right)$ which generates $M$ as a $\Ulie(\gfrak)$-module.  For every $n\in\amsbb{Z}_{>0}$, let $\gfrak'_n$ denote the subalgebra $\hfrak+\gfrak_n$ and set $\bfrak_n'\subseteq\gfrak_n'$ to be the Borel subalgebra $\hfrak+\bfrak_n$ of $\gfrak_n'$.  Then, the $\Ulie\left(\gfrak'_n\right)$-module  $M_n$ is given by $
				M_n:=\sum_{i=0}^n\,\Ulie\left(\gfrak'_n\right)\cdot v_i$.
	
	Note that the finite-dimensional theory carries trivially over to $\gfrak_n'$, and we use the notation $\bggO^{\gfrak'_n}_{\bfrak'_n}$ for the category $\bggO$ of $\gfrak'_n$ with respect to the Borel subalgebra $\bfrak_n'$.  Denote by $\bggO^{\gfrak'_n}_{\bfrak'_n}\left\llbracket\lambda\right\rrbracket$ the block of $\bggO^{\gfrak'_n}_{\bfrak'_n}$ containing $\Llie_n\left(\lambda\right):=\Llie\left(\lambda;\gfrak'_n,\bfrak'_n,\hfrak\right)$.  Note that we have the direct sum decomposition $M_n=X_n\oplus Y_n$, where $X_n\in \bggO^{\gfrak_n}_{\bfrak_n}\left\llbracket\lambda\right\rrbracket$ and all composition factors of $Y_n$ are not in $ \bggO^{\gfrak'_n}_{\bfrak'_n}\left\llbracket\lambda\right\rrbracket$.  The submodules $X_n$ and $Y_n$ are unique.  Furthermore, if $N_n$ is an indecomposable submodule of $M_n$, then $N_n$ must lie entirely in $X_n$ or in $Y_n$.
	
	 Define 
	\begin{align} X:=\big\{x\in M\,\boldsymbol{|}\,x\in X_n\textup{ for all sufficiently large }n\big\}\end{align}
	and
	\begin{align}
		Y:=\textup{span}_\amsbb{K}\big\{y\in M\,\boldsymbol{|}\,y\in Y_n\textup{ for infinitely many }n\big\}\,.
	\end{align}
	Then, it is evident that $X$ and $Y$ are $\gfrak$-submodules of $M$.  We shall  prove  that $X+Y=M$ and that $X\cap Y=0$.  
	
	First, let the $\Ulie\left(\gfrak'_n\right)$-module $X_n'$ be an indecomposable direct summand of $X_n$.  Fix $\tilde{n}\geq n$.   Note that we have either $X'_n\subseteq X_{\tilde{n}}$ or $X'_n\subseteq Y_{\tilde{n}}$.   Likewise, if $Y_n'$ is an indecomposable direct summand of $Y_n$, then either $Y'_n\subseteq Y_{\tilde{n}}$ or $Y'_n\subseteq X_{\tilde{n}}$.	
	
	To justify the statement in the paragraph above, consider the $\Ulie\left(\gfrak'_{\tilde{n}}\right)$-module 
	$
		\tilde{X}_{\tilde{n}}:=\Ulie\left(\gfrak'_{\tilde{n}}\right)\cdot X'_n$.
	Let $\Xi\subseteq\hfrak^*$ denote the set of weights of $X'_n$.  By (\ref{eq:equaldim}), $\Xi$ is also an indecomposable weight set of $\tilde{X}_{\tilde{n}}$.  Therefore, $\tilde{X}_{\tilde{n}}$ can be decomposed as a direct sum $\tilde{X}_{\tilde{n}}^1\oplus\tilde{X}_{\tilde{n}}^2\oplus\ldots\oplus\tilde{X}^l_{\tilde{n}}$, where each $\tilde{X}^i_{\tilde{n}}$ is an indecomposable $\Ulie\left(\gfrak'_{\tilde{n}}\right)$-module, but as $\Xi$ is an indecomposable weight set of $\tilde{X}_{\tilde{n}}$, we must have $\Xi\subseteq \textup{supp}\left(\tilde{X}^i_{\tilde{n}}\right)$ for some $i$.  However, this means $X'_n\subseteq \tilde{X}^i_{\tilde{n}}$.  Now, being indecomposable, $\tilde{X}^i_{\tilde{n}}$ must lie entirely either in $X_{\tilde{n}}$ or in $Y_{\tilde{n}}$.  Ergo, $X'_n$ is a subspace of $X_{\tilde{n}}$ or $Y_{\tilde{n}}$ for every $\tilde{n}\geq n$.

	The paragraph above proves that $X\cap X_n$ is given by a direct sum of some indecomposable direct summands of $X_n$.  Indeed, for a fixed direct summand $X'_n$ of $X_n$, we have only two possible scenarios: either $X'_n$ lies in $X_{\tilde{n}}$ for all sufficiently large $\tilde{n}\geq n$, or $X'_n$ lies in $Y_{\tilde{n}}$ for infinitely many $\tilde{n}$.  In the former case, $X'_n\subseteq X$, whereas, in the latter case, $X'_n\cap X=0$ and $X'_n\subseteq Y$.  In other words, $A_n:=X\cap X_n$ is a direct summand of $M_n$.  Write $B_n:=Y_n\oplus Z_n$, where $Z_n$ is the direct sum of indecomposable direct summands $X'_n$ of $X_n$ which intersect $X$ trivially. 
	
	Next, we fix $\xi\in\textup{supp}(M)$.  We shall verify that $M^\xi=X^\xi+Y^\xi$.  For a given $v\in M^\xi$, $v=a_n+b_n$ for some $a_n\in \left(A_n\right)^{\xi}$ and $b_n\in \left(B_n\right)^{\xi}$.  Suppose that $n_0$ is a positive integer such that $\left(M_n\right)^{\xi}=M^\xi$ for all $n\geq n_0$.    We claim that there exists a positive integer $n_1\geq n_0$ such that $a_{n_1}=a_{n_1+1}=a_{n_1+2}=\ldots$.   This claim follows from the observation that \begin{align}A_n\subseteq A_{n+1}\subseteq A_{n+2}\subseteq \ldots\end{align} for all $n\geq n_0$.  The finite-dimensionality assumption implies that 
	\begin{align}
		\left(A_{n_1}\right)^{\xi}=\left(A_{n_1+1}\right)^{\xi}=\left(A_{n_1+2}\right)^{\xi}=\ldots\end{align} for some $n_1\geq n_0$.  Furthermore, we note that 
	\begin{align}
		B_n\supseteq B_{n+1}\supseteq B_{n+2}\supseteq\ldots\,;
	\end{align} consequently, the finite-dimensionality assumption yields \begin{align}\left(B_{n_1}\right)^{\xi}=\left(B_{n_1+1}\right)^{\xi}=\left(B_{n_1+2}\right)^{\xi}=\ldots\,.\end{align}    The claim follows immediately. 
	
	We write $a$ for the common value $a_{n_1}=a_{n_1+1}=a_{n_1+2}=\ldots$.  Set $b:=v-a$.  We shall now justify that $b$ is an element of $Y$.  Recall that $B_n=Y_n\oplus Z_n$, where $Z_n$ is the direct sum of indecomposable direct summands of $X_n$ that intersect $X$ trivially.  For $n\geq n_1$, we can write 
	\begin{align}
		b=\sum_{i=1}^{k_n}\,y_n^i+\sum_{j=1}^{l_n}\,z_n^j\,,
	\end{align}
	where $y_n^i$ and $z_n^j$ are nonzero elements of indecomposable direct summands of $Y_n$ and $Z_n$.  We shall now prove that each $y_n^i$ and each $z_n^j$ belong in $Y$.  For $\tilde{n}\geq n$, note that each $y_n^i$ lies either in $X_{\tilde{n}}$ or in $Y_{\tilde{n}}$.  If the former scenario occurs for all sufficiently large $\tilde{n}\geq n$, then $y_n^i\in X$, but this immediately implies $y_n^i=0$, which is a contradiction.  Ergo, the latter scenario occurs for infinitely many values $\tilde{n}\geq n$, whence $y_n^i \in Y$.  The same argument applies to each $z_n^j$.  Thus, we conclude that $y_n^i,z_n^j\in Y$ for every $i=1,2,\ldots,k_n$ and $j=1,2,\ldots,l_n$ with $n\geq n_1$.  Thence, $b\in Y$.  This proves that $M^\xi=X^\xi+Y^\xi$, leading to $M=X+Y$.
	
	Now, we shall check that $X\cap Y=0$.  Let $y_1,y_2,\ldots,y_k $ be linearly independent elements of $Y$ such that $x:=y_1+y_2+\ldots+y_k$ is in $X$ and that, for each $i=1,2,\ldots,k$, there are infinitely many positive integers $n$ for which $y_j\in Y_n$.  We may assume that there exists $\xi\in\textup{supp}(M)$ with $y_i\in M^\xi$ for every $i=1,2,\ldots,k$.   Additionally, there exists a positive integer $m$ such that $x\in M_m$ and that $x\in X_n$ for every integer $n\geq m$.  
	
	Let $\bar{n}(j) \geq m$ be a positive integer such that $y_j\in Y_{\bar{n}(j)}$.  We decompose $y_j$ as
	\begin{align}
		y_j = y_j^1+y_j^2+\ldots+y_j^{r_j}\,,
	\end{align}
	where each $y_j^i$ is nonzero and in an indecomposable direct summand of $Y_{\bar{n}(j)}$.
	Pick an arbitrary $n\geq\bar{n}(j)$.  We note that each $y_j^i$ must lie in $X_n$ or in $Y_n$.  However, as $x\in X_n$, we conclude that $y_j^i$ is in  $X_n$, whence $y^i_j\in X_n$ for every $n\geq \bar{n}(j) $.  As a result, $y_j=\sum_{i=1}^{r_j}\,y_j^i$ is in $X$, which means $y_j=0$, and a contradiction is reached.

	Finally, we have the following four equalities: $\big[X:\Llie(\lambda)\big]=\big[M:\Llie(\lambda)\big]>0$, $\big[X:\Llie(\mu)\big]=0$, $\big[Y:\Llie(\lambda)\big]=0$, and $\big[Y:\Llie(\mu)\big]=\big[M:\Llie(\mu)\big]>0$.  That is, $M=X\oplus Y$ with $X\neq 0$ and $Y\neq 0$ .  This contradicts the assumption that $M$ is indecomposable.
		\end{pf} 
	\end{prop}
	

	\begin{prop}
		A block of $\bbggO$ containing $\Llie(\lambda)$ contains $\bbggO\llbracket\lambda\rrbracket$ as a subcategory.
	
		\begin{pf}
			Using the indecomposability of the Verma modules, we conclude that the block containing $\Llie(\lambda)$ must have $\bar\bggO\llbracket\lambda\rrbracket$ as a subcategory.  In other words, let $\mu,\nu\in \llbracket\lambda\rrbracket$.  Let $n\in\amsbb{N}$ be sufficiently large that $\mu=w \cdot \lambda$ for some $w\in W_n\left[\lambda_n\right]$.   (Here, $\xi_n$ denotes $\xi|_{\hfrak_n}$ for all $\xi\in\hfrak^*$.)
			
			From the finite-dimensional theory (see~\cite{bggo}), $W_n\left[\lambda_n\right]\cdot\lambda$ has a unique maximal element $\upsilon$ (with respect to the order $\preceq$ given by $\bfrak$).  Then, the Verma module $\Mlie(\upsilon)$ has $\Mlie(\mu)$ and $\Mlie(\lambda)$ as submodules due to the BGG Theorem (Theorem~\ref{thm:stronglinkagegen}).  Therefore, we have nonzero homomorphisms $\Mlie(\mu)\to \Mlie(\upsilon)$ and $\Mlie(\lambda)\to \Mlie(\upsilon)$.   Thus, the indecomposable modules $\Mlie(\mu)$ and $\Mlie(\lambda)$ are in the same block.  Furthermore, with  nontrivial homomorphisms $\Mlie(\mu)\to \Llie(\mu)$ and $\Mlie(\lambda)\to\Llie(\lambda)$, we conclude that $\Llie(\mu)$ and $\Llie(\lambda)$ are in the same block.
			
			Now, suppose that $M\in\bbggO$ is indecomposable with $\Llie(\mu)$ as a composition factor, where $\mu\in \llbracket\lambda\rrbracket$.  By Theorem~\ref{thm:gencomp} and Proposition~\ref{prop:blocks}, we see that $M$ has a submodule $N$ such that $M/N\cong \Llie(\nu)$ for some $\nu \in \llbracket\lambda\rrbracket$.  Thus, the nonzero homomorphism $M\to M/N$ establishes that $M$ is in the same block as $\Llie(\nu)$, which is also in the same block as $\Llie(\lambda)$.  Thus, every indecomposable object $M$ whose composition factors are of the form $\Llie(\mu)$ with $\mu\in \llbracket\lambda\rrbracket$ is in the same block as $\Llie(\lambda)$.  The proposition follows immediately.
		\end{pf}
	\end{prop}

	\begin{thm}
		Let $\Omega$ denote the set of integral Weyl dot-orbits.  Then, the full abelian subcategories $\bar\bggO\llbracket\lambda\rrbracket$, where $\llbracket\lambda\rrbracket\in\Omega$, are the blocks of $\bar\bggO$, and
		$
		\bar	\bggO =\bigoplus_{\llbracket\lambda\rrbracket\in \Omega}\,\bar\bggO\llbracket\lambda\rrbracket$.
		\label{thm:blocks}
		
		\begin{pf}
			From the proposition above, we know that each block of $\bbggO$ contains $\bbggO\llbracket\lambda\rrbracket$ for some weight $\lambda\in\hfrak^*$.  We shall prove that the block containing $\bbggO\llbracket\lambda\rrbracket$ must then coincide with $\bbggO\llbracket\lambda\rrbracket$.  If the block contains an indecomposable module $M$ which is not in $\bbggO\llbracket\lambda\rrbracket$, then there exists a finite sequence $M=M_0,M_1,M_2,\ldots,M_k$ of indecomposable modules, all of which belong in this block, with the properties that $M_k\in\bbggO\llbracket\lambda\rrbracket$, $M_{k-1}\notin\bbggO\llbracket\lambda\rrbracket$, and for each $i=0,1,2,\ldots,k-1$, either $\Hom_{\bbggO}\left(M_i,M_{i+1}\right)\neq 0$ or $\Hom_{\bbggO}\left(M_{i+1},M_i\right)\neq 0$.
			
			If $\Hom_{\bbggO}\left(M_{k-1},M_{k}\right)\neq 0$ or $\Hom_{\bbggO}\left(M_{k},M_{k-1}\right)\neq 0$, then $M_{k-1}$ has a composition factor $\Llie(\mu)$ (which is also a composition factor of $M_k$) for some $\mu\in \llbracket\lambda\rrbracket$, which then means that $M_{k-1}\in\bbggO\llbracket\lambda\rrbracket$ by Proposition~\ref{prop:blocks}.  This contradicts the assumption that $M_{k-1}\notin \bbggO\llbracket\lambda\rrbracket$.  Therefore, the blocks of $\bbggO$ are precisely $\bbggO\llbracket\lambda\rrbracket$.
			
			To complete the proof, let now $M$ be an arbitrary object in $\bbggO$.  By Theorem~\ref{thm:directsumdecomposition}, $M$ has a direct sum decomposition with indecomposable summands.  Write $M\llbracket\lambda\rrbracket$ for the (direct) sum of the direct summands of $M$ that belong to $\bbggO\llbracket\lambda\rrbracket$.  Then, we can clearly see that
		$
			M=\bigoplus_{\llbracket\lambda\rrbracket\in \Omega}\,M\llbracket\lambda\rrbracket$.
		Note that this direct sum may be an uncountable direct sum.
		\end{pf}	
	\end{thm}

\subsection{Kazhdan-Lusztig Multiplicities}

Let $\gfrak_n':=\gfrak_n+\hfrak$ and $\bfrak'_n:=\bfrak_n+\hfrak$.  Note that the Weyl group of $\gfrak_n'$ is still the Weyl group $W_n$ of $\gfrak_n$.  For each $\xi\in\hfrak^*$, write $\Mlie_n(\xi)$ and $\Llie_n(\xi)$ for the Verma module $\Mlie\left(\xi;\gfrak'_n,\bfrak_n',\hfrak\right)$ and the simple module $\Llie\left(\xi;\gfrak_n',\bfrak_n',\hfrak\right)$, respectively.

\begin{define}
	A weight $\lambda\in\hfrak^*$ is \emph{regular} if $(\lambda+\rho)\left(h_\alpha\right)\neq 0$ for every root $\alpha$.
\end{define}

Fix a regular integral weight $\lambda$.  Take $\mu\in W\cdot \lambda$.  For each $n\in\amsbb{Z}_{> 0}$, write $\nu_n$ for the antidominant weight in $\hfrak^*$ that is strongly linked to $\lambda$ with respect to $\bfrak_n'$.  In addition, there exist elements $x_n$ and $y_n$ of $W_n$ such that 
$ x_n^{-1}\cdot\lambda=  \nu_n$ and $y_n^{-1}\cdot \mu=\nu_n$
From the finite-dimensional theory, we have
\begin{align}
	\big[\Mlie_n\left(\lambda\right):\Llie_n\left(\mu\right)\big]=P^{W_n}_{w_n^0x_n,w_n^0y_n}(1)\,,
\end{align}
where $w_n^0\in W_n$ is the longest element of $W_n$.  Combining this result with the observation that $\big[\Mlie_n\left(\lambda\right):\Llie_n\left(\mu\right)\big]=\big[\Mlie(\lambda),\Llie(\mu)\big]$ for sufficiently large values of $n$, we obtain the proposition below.

\begin{prop}
	There exists a positive integer $n_0$ such that, for all $n\geq n_0$,
	\begin{align}
		\big[\Mlie\left(\lambda\right):\Llie\left(\mu\right)\big]=P^{W_n}_{w_n^0x_n,w_n^0y_n}(1)\,.
	\end{align}
\end{prop}

 Fix $x,y\in W$ and set $m(x,y)$ to be the smallest positive integer $m$ such that $x,y\in W_m$.   From Chapter 1.10 of \cite{Coxeter}, we see that
 $
 	P_{x,y}^{W_{m(x,y)}}(q)=P_{x,y}^{W_n}(q)=P_{x,y}^W(q)$.
 This result gives rise to the following proposition.
 
 \begin{prop}
 	For every $x\in W$ and for each regular antidominant weight $\lambda$, 
 	\begin{align}
 		\big[\Mlie(x\cdot \lambda)\big]=\sum_{y\preceq x}\,P^{W_{m(x,y)}}_{w^0_{m(x,y)}x,w^0_{m(x,y)}y}(1)\,\big[\Llie(y\cdot \lambda)\big]\,,
 	\end{align}
 	or equivalently
 	\begin{align}
 		\big[\Llie(x\cdot \lambda)\big]=\sum_{y\preceq x}\,(-1)^{\ell(x)-\ell(y)}\,P^{W_{m(x,y)}}_{x,y}(1)\,\big[\Mlie(y\cdot\lambda)\big]\,.
 	\end{align}
 	(Note that the two equations above are equalities in the Grothendieck group of $\bar\bggO$.)
 \end{prop}
 
 \subsection{$\Hom$ and $\Ext^\bullet$ Functors}
 
 	Unless otherwise specified, $\Ext$ denotes $\Ext_{\bar\bggO}$.  Similarly, $\Hom$ denotes $\Hom_{\bar\bggO}$.
 	
 	\begin{prop}
 		Let $\lambda,\mu\in\hfrak^*$.
 		\label{prop:extintro}
 		
 		\begin{enumerate}[(a)]
 			\item If a $\gfrak$-module $M$ is such that, for all $\upsilon\in \textup{supp}(M)$, we have $\lambda \not\prec \upsilon$, then $\Ext^1\big(\Mlie(\lambda),M\big)=0$.  In particular,
 			\begin{align}
 				\Ext^1\big(\Mlie(\lambda),\Llie(\lambda)\big)=0\textup{ and }\Ext^1\big(\Mlie(\lambda),\Mlie(\lambda)\big)=0\,.
 			\end{align}
 			\item If $\mu\preceq \lambda$, then $\Ext^1\big(\Mlie(\lambda),\Llie(\mu)\big)=0$.
 			\item If $\mu\prec\lambda$ and $\Nlie(\lambda)$ is the maximal proper submodule of $\Mlie(\lambda)$, then
 			\begin{align}
 				\Ext^1\big(\Llie(\lambda),\Llie(\mu)\big)\cong \Hom\big(\Nlie(\lambda),\Llie(\mu)\big)\,.
 				\label{eq:extandN}
 			\end{align}
 			\item $\Ext^1\big(\Llie(\lambda),\Llie(\lambda)\big)=0$.
 		\end{enumerate}
 		
 		\begin{pf}$\phantom{a}$
 		
 			\begin{enumerate}[(a)]
 				\item
 				Given an extension $0 \to M \overset{i}\longrightarrow E \overset{p}{\longrightarrow} \Mlie(\lambda) \to 0$ in $\bar\bggO$, let $e\in E$ be such that $p(e)$ be a highest-weight vector of $\Mlie(\lambda)$.  Due to the hypothesis, the submodule $V$ of $E$ generated by $e$ is a highest-weight module with highest weight $\lambda$.  Since $V$ is mapped surjectively by $p$ onto $\Mlie(\lambda)$, we conclude that $p$ induces an isomorphism $V \cong \Mlie(\lambda)$, whence the exact sequence splits.
 				\item This is an immediate consequence of (a).
 				\item Starting with the short exact sequence 
				$0\to \Nlie(\lambda) \to \Mlie(\lambda) \to \Llie(\lambda)\to 0$, we get the following long exact sequence of $\Ext$-groups:
 				\begin{align}
 					\ldots \to \Hom\big(\Mlie(\lambda),\Llie(\mu)\big) &\to \Hom\big(\Nlie(\lambda),\Llie(\mu)\big)
 					\nonumber\\
 					&\to \Ext^1\big(\Llie(\lambda),\Llie(\mu)\big) \to \Ext^1\big(\Mlie(\lambda),\Llie(\mu)\big)\to\ldots\,.
 				\end{align}
 				By (b), $\Ext^1\big(\Mlie(\lambda),\Llie(\mu)\big)=0$.  Furthermore, it is clear that $\Hom\big(\Mlie(\lambda),\Llie(\mu)\big)=0$.  Therefore, we have the isomorphism (\ref{eq:extandN}).
 				\item Replace $\mu$ by $\lambda$ in the proof of (c).  We note that $\Hom\big(\Nlie(\lambda),\Llie(\lambda)\big)=0$.  By (b), we have $\Ext^1\big(\Mlie(\lambda),\Llie(\mu)\big)=0$.  Thus, $\Ext^1\big(\Llie(\lambda),\Llie(\lambda)\big)=0$ as well.\qedhere
 			\end{enumerate}
 		\end{pf}
 	\end{prop}
 	
 	\begin{prop}
 		Let $\lambda,\mu\in\hfrak^*$.
 		\begin{enumerate}[(a)]
 			\item For every $M,N\in\bar\bggO$ and $k\in\amsbb{Z}_{\geq 0}$, we have $\Ext^k(M,N)\cong \Ext^k\big(N^\vee,M^\vee\big)$.
 			\item The image of any homomorphism $\Mlie(\mu)\to\Vlie(\lambda)$ is a submodule of $\Llie(\lambda)\subseteq \Vlie(\lambda)$.  This means 
			\begin{align}
				\dim_\amsbb{K}\Big(\Hom\big(\Mlie(\mu),\Vlie(\lambda)\big)\Big)=\left\{
					\begin{array}{ll}
						1\,,	&\textup{if }\mu=\lambda\,,
						\\
						0\,,	&\textup{if }\mu\neq \lambda\,.
					\end{array}
				\right.
			\end{align}
 			\item $\Ext^1\big(\Mlie(\mu),\Vlie(\lambda)\big)=0$ for all $\lambda$ and $\mu$.
 		\end{enumerate}
 		\label{prop:extverma}
 		
 		\begin{pf}$\phantom{a}$
 			
 			\begin{enumerate}[(a)]
 				\item This part is trivial due to the fact that duality is an antiequivalence of the category $\bbggO$ with itself.
 				\item Let $M$ be the image of a nonzero homomorphism $\Mlie(\mu)\to\Vlie(\lambda)$.  Then, $M$ is a highest-weight submodule of $\Vlie(\lambda)$ with highest weight $\mu$.  Since $\Llie(\lambda)$ is contained in every nonzero submodule of $\Vlie(\lambda)$, we see that $\Llie(\lambda)\subseteq M$, so $\mu\preceq \lambda$.  However, the composition factors of $\Vlie(\lambda)$ are the same as those of $\Mlie(\lambda)$, which are simple modules with highest weight less than or equal to $\lambda$.  This means $\mu\succeq \lambda$.  Consequently, $\mu=\lambda$ must hold, whence $M=\Llie(\lambda)$.
 				\item   If $\lambda\not\prec \mu$, then $M:=\Vlie(\mu)$ satisfies the hypothesis of Proposition~\ref{prop:extintro}(a).  Therefore, 
 				\begin{align}\Ext^1\big(\Mlie(\lambda),\Vlie(\mu)\big)=0\,.\end{align}  By (a), we have $
 						\Ext^1\big(\Mlie(\mu),\Vlie(\lambda)\big)\cong\Ext^1\big(\Mlie(\lambda),\Vlie(\mu)\big)=0$.
					If $\lambda \preceq \mu$, then $M:=\Vlie(\lambda)$ satisfies the hypothesis of Proposition~\ref{prop:extintro}(a), with $\mu$ replacing $\lambda$ in that proposition.  The same conclusion follows.
 			\end{enumerate}
 		\end{pf}
 	\end{prop}

%
%

\section{Truncated Category $\bggO$}
\label{ch:trunc}


\subsection{Truncation}
	
	As before, $M^\vee$ and $f^\vee$ denote the duals  in $\bar\bggO$ of an object $M$ and a homomorphism $f$, respectively.   We shall now define a truncation method of the category $\bbggO$ using an idea from~\cite{RCW}.

\begin{define}
	 For $\lambda\in\hfrak^*$, we write $\bar\bggO^{\preceq \lambda}$ for the full subcategory of $\bar\bggO$ consisting of all modules $M\in\bar\bggO$ whose weights are less than or equal to $\lambda$ with respect to the partial order $\preceq$ on $\hfrak^*$.
\end{define}

\begin{prop}
	For each $\lambda\in\hfrak^*$, let $\tcyr_{\preceq \lambda}:\bar\bggO\to\bar\bggO^{\preceq\lambda}$ be defined as 
	\begin{align}
		\tcyr_{\preceq\lambda}M\defeq\sum_{\substack{{N\subseteq M}\\{N\in\bar\bggO^{\preceq\lambda}}}}\,N
	\textup{ and }
		\tcyr_{\preceq\lambda}f \defeq f|_{\tcyr_{\preceq\lambda}M}
	\end{align}
	for all $M\in\bar\bggO$ and for all homomorphisms $f:M\to L$ of objects in $\bar\bggO$.
	Then, $\tcyr_{\preceq\lambda}$ is a left-exact (covariant) functor.  We shall call $\tcyr_{\preceq \lambda}$ the \emph{truncation functor} (with the upper bound $\lambda$).
\end{prop}

\begin{cor}
	For each $\lambda\in\hfrak^*$, let $\tcyr^\vee_{\preceq \lambda}:\bar\bggO\to\bar\bggO^{\preceq\lambda}$ be defined as 
	\begin{align}
		\tcyr^\vee_{\preceq\lambda}M\defeq \Big(\tcyr_{\preceq\lambda}\left(M^\vee\right)\Big)^\vee
\textup{ and }
		\tcyr^\vee_{\preceq\lambda}f \defeq \Big(\tcyr_{\preceq\lambda}\left(f^\vee\right)\Big)^\vee
	\end{align}
	for all $M\in\bar\bggO$ and for all homomorphisms $f:M\to L$ of objects in $\bar\bggO$.
	Then, $\tcyr^\vee_{\preceq\lambda}$ is a right-exact (covariant) functor.  We shall call $\tcyr^\vee_{\preceq \lambda}$ the \emph{dual truncation functor} (with the upper bound $\lambda$).
	\label{cor:modtrunc}
\end{cor}

\begin{prop}
	Let $\lambda\in\hfrak^*$.  If $I$ is an injective object in $\bar\bggO$, then $\tcyr_{\preceq\lambda}I$ is injective in $\bar\bggO^{\preceq \lambda}$.  If $P$ is a projective object in $\bar\bggO$, then $\tcyr^\vee_{\preceq\lambda}P$ is projective in $\bar\bggO^{\preceq \lambda}$. 
	
	\begin{pf}
		Let $I$ be an injective object in $\bar\bggO$ and $0\to \tcyr_{\preceq\lambda}I \to M\to N\to 0$ an exact sequence of objects in $\bar\bggO^{\preceq\lambda}$.  We have the injection $\tcyr_{\preceq\lambda}I \overset{\subseteq}{\longrightarrow} I$.  Because $I$ is an injective object in $\bar\bggO$ and $0\to \tcyr_{\preceq\lambda}I \to M\to N\to 0$ is also an exact sequence of objects in $\bar\bggO$, we conclude that there exists a homomorphism $\phi:M\to I$ such that the diagram below is commutative:
		
		\begin{equation}
		\begin{tikzcd}
		0\arrow{r}& \tcyr_{\preceq \lambda}I \arrow{r} \arrow{d}{\rotatebox{270}{$\subseteq$} }& M \arrow{r}\arrow{ld}{\phi} &N \arrow{r} & 0\\
		&I\,.&&&
		\end{tikzcd}
		\end{equation}

		Since the image of $M$ under $\phi$ must be in $\bar\bggO^{\preceq\lambda}$, we see that $\im(\phi)\subseteq \tcyr_{\preceq\lambda}I$.  Thence, we indeed have a commutative diagram
		
		\begin{equation}
		\begin{tikzcd}
		0\arrow{r}& \tcyr_{\preceq \lambda}I \arrow{r} \arrow{d}{\rotatebox{270}{$=$} }& M \arrow{r}\arrow{ld}{\phi} &N \arrow{r} & 0\\
		&\tcyr_{\preceq \lambda}I \,.&&&
		\end{tikzcd}
		\end{equation}
		Hence, the exact sequence  $0\to \tcyr_{\preceq\lambda}I \to M\to N\to 0$ splits.  Thus, $\tcyr_{\preceq\lambda}I$ is injective.  
		
		For the second part of the proposition, we employ the duality from Corollary~\ref{cor:modtrunc}.  The proof is now complete.
	\end{pf}
\end{prop}


\subsection{Injective Objects}
\label{sec:injectiveobjects}

As in the proof of Proposition~\ref{prop:blocks}, $\gfrak'_n$ and $\bfrak'_n$ denote $\hfrak+\gfrak_n$ and $\hfrak+\bfrak_n$, respectively.

\begin{prop}
		Fix $n\in\amsbb{Z}_{> 0}$.  Suppose that $I_{n+1}$ is an injective object in $\bggO^{\gfrak'_{n+1}}_{\bfrak'_{n+1}}$.  Then, the restriction $I_n:=\Res^{\gfrak'_{n+1}}_{\gfrak'_n}\left(I_{n+1}\right)$ is an injective object in $\bggO^{\gfrak'_n}_{\bfrak'_n}$.		
		\label{prop:injrestricted}
		
		\begin{pf}
			Let $0\to M_n\overset{\varphi_n}{\longrightarrow} N_n$ be an exact sequence in $\bggO^{\gfrak'_n}_{\bfrak'_n}$ along with a homomorphism $f_n:M_n\to I_{n}$ of $\gfrak'_n$-modules.  Now, let $\pfrak'_{n+1}$ denote the parabolic subalgebra
			$
				\pfrak'_{n+1}=\gfrak'_n+\bfrak'_{n+1}$.
			Take $\left\{x_{\pm\alpha}\,\boldsymbol{|}\,\alpha\in\Delta^+\right\}\cup\left\{h_\beta\,\boldsymbol{|}\,\beta\in\Sigma^+\right\}$ for a Chevalley basis of $\gfrak$.  We equip each object $L$ in $\bggO^{\gfrak'_n}_{\bfrak'_n}$ with a $\pfrak'_{n+1}$-module structure by requiring that,  for each $\bfrak'_{n+1}$-positive root $\alpha$ of $\gfrak'_{n+1}$ which is not a root of $\gfrak'_n$, 
			$
				x_\alpha\cdot L=0$.

			Note that $\Ulie\left(\gfrak'_{n+1}\right)$ is a free (whence flat) $\Ulie\left(\pfrak'_{n+1}\right)$-module due to the PBW Theorem.  Hence, the parabolic induction functor $\Ulie\left(\gfrak'_{n+1}\right)\underset{\Ulie\left(\pfrak'_{n+1}\right)}{\otimes}\boldsymbol{\_}$ is exact, that is we have an exact sequence of $\gfrak'_{n+1}$-modules 
			\begin{align}
				0\to \Ulie\left(\gfrak'_{n+1}\right)\underset{\Ulie\left(\pfrak'_{n+1}\right)}{\otimes}M_n\overset{\varphi_{n+1}}{\longrightarrow} \Ulie\left(\gfrak'_{n+1}\right)\underset{\Ulie\left(\pfrak'_{n+1}\right)}{\otimes}N_n\,,
				\end{align}
			where the $\gfrak_{n+1}'$-module homomorphism $\varphi_{n+1}$ is given by $\varphi_{n+1}:=\textup{id}_{\Ulie\left(\gfrak_{n+1}'\right)}\otimes \varphi_n$.  Then, we define the $\gfrak'_{n+1}$-module homomorphism $f_{n+1}:\Ulie\left(\gfrak_{n+1}'\right)\underset{\Ulie\left(\pfrak'_{n+1}\right)}{\otimes}M_n\to I_{n+1}$ by setting
			$
				f_{n+1}(x\otimes v):=x\cdot f_n(v)
			$
			for all $x\in \Ulie\left(\gfrak'_{n+1}\right)$ and $v\in M_n$.  
			
			Since $\Ulie\left(\gfrak'_{n+1}\right)\underset{\Ulie\left(\pfrak'_{n+1}\right)}{\otimes}M_n$ and $\Ulie\left(\gfrak'_{n+1}\right)\underset{\Ulie\left(\pfrak'_{n+1}\right)}{\otimes}N_n$ are objects in the category $\bggO^{\gfrak'_{n+1}}_{\bfrak'_{n+1}}$, we conclude by injectivity of the module $I_{n+1}$ in $\bggO^{\gfrak'_{n+1}}_{\bfrak'_{n+1}}$ that there exists a homomorphism of $\gfrak'_{n+1}$-modules $\psi_{n+1}:\Ulie\left(\gfrak'_{n+1}\right)\underset{\Ulie\left(\pfrak'_{n+1}\right)}{\otimes}N_n\to I_{n+1}$ such that 
			$
				f_{n+1}=\psi_{n+1}\circ\varphi_{n+1}
			$.
			We then define $\psi_n:N_n\to I_{n+1}$ by setting 
			$
				\psi_n(u):=\psi_{n+1}\left(1_{\Ulie\left(\gfrak'_{n+1}\right)}\otimes u\right)
			$
			for every $u\in N_n$.  It is easy to see that $f_n=\psi_n\circ \varphi_n$ and that $I_n=\Res^{\gfrak'_{n+1}}_{\gfrak'_n}\left(I_{n+1}\right)$ is an object in $\bggO^{\gfrak'_n}_{\bfrak_n'}$, whence $I_{n}$ is injective in $\bggO^{\gfrak'_n}_{\bfrak'_n}$.
		\end{pf}
	\end{prop}

\begin{define}
Let $\lambda\in\hfrak^*$.  Then, we say that $\lambda$ is \emph{dominant} if $(\lambda+\rho)\left(h_\alpha\right)\notin\amsbb{Z}_{<0}$ for all positive roots $\alpha$ of $\gfrak$.  We say that $\lambda$ is \emph{almost dominant} if $(\lambda+\rho)\left(h_\alpha\right)\in\amsbb{Z}_{<0}$ for at most finitely many positive roots $\alpha$ of $\gfrak$.
\label{def:dominantweights}
\end{define}

\begin{thm}
	Let $I_n\in\bggO^{\gfrak'_n}_{\bfrak'_n}$ for each $n\in\amsbb{Z}_{> 0}$ be an injective object.  Suppose that we have an embedding $I_n\overset{\psi_n}{\longrightarrow}I_{n+1}$ for every $n$.  If the direct limit $I=\lim\limits_{\underset{n}{\boldsymbol{\longrightarrow}}}\,I_n$ is an object of $\bar\bggO$, then $I$ is injective in $\bar\bggO$.
	\label{thm:injdirlim}
	
	\begin{pf}
		Let an injective homomorphism $M\overset{\varphi}{\longrightarrow}N$ and a homomorphism $f:M\to I$ be given.  Without loss of generality, we can assume that all weights of $M$, $N$, and $I$ lie within $\lambda+\Lambda$ for some $\lambda\in\hfrak^*$.  In particular, we can assume that the modules $M$, $N$, and $I$ are generated by countably many weight vectors.
		
		We suppose that $N$ is generated by weight vectors $u_1,u_2,u_3,\ldots$, with the corresponding weights $\mu_1,\mu_2,\mu_3,\ldots$.  Let $N_n$ denote the $\gfrak_n$-submodule $
			N_n:=\sum_{i=1}^n\,\Ulie\left(\gfrak'_n\right)\cdot u_i$.
		Then, we define $M_n$ as $\varphi^{-1}\left(N_n\right)$.  Note that both $M_n$ and $N_n$ are objects in $\bggO^{\gfrak'_n}_{\bfrak'_n}$.  By Proposition~\ref{prop:injrestricted}, for every $m\geq n$, $\textup{Res}^{\gfrak'_m}_{\gfrak'_n}\left(I_m\right)$ is an injective module in $\bggO^{\gfrak'_n}_{\bfrak'_n}$.  Therefore, as $I$ has finite-dimensional weight spaces, we can assume without loss of generality that each $I_n$ satisfies the condition $
			\dim_\amsbb{K}\big(\left(I_n\right)^{\mu_i}\big)=\dim_\amsbb{K}\left(I^{\mu_i}\right)$
		for every $i=1,2,\ldots,n$.
		
		Let $\varphi_n:=\varphi|_{M_n}$ and $f_n:=f|_{M_n}$.  From the definitions above, we have the diagram of objects of $\bggO^{\gfrak'_n}_{\bfrak'_n}$
		\begin{equation}
		\begin{tikzcd}
		0\arrow{r}& M_n \arrow{r}{\varphi_n} \arrow{d}{f_n} & N_n\\
		&I_n \,.&
		\end{tikzcd}
		\end{equation}
		Because the object $I_n$ is injective, there exists a $\gfrak'_n$-module homomorphism $F_n:N_n\to I_n$ such that $F_n\circ \varphi_n=f_n$.
		
		Write $\mathcal{F}^1_n$ for the set of all maps $F_n:N_n\to I_n$ such that $
			F_n\circ \varphi_n=f_n$.
		Take $V^1_n$ to be the $\amsbb{K}$-span of all  vectors of the form $\big(1,F_n\left(u_1\right)\big) \in \amsbb{K}\times I^{\mu_1}$, where $F_n\in \mathcal{F}^1_n$.    Then, $V_n^1$ is a finite-dimensional vector space for every $n$.  Furthermore, we have $
			V_1^1\supseteq V_2^1\supseteq V_3^1\supseteq \ldots$.
		Hence, the inclusion sequence above stabilizes at some $n_1\in\amsbb{Z}_{>0}$.  That is,
		$
			V^1:=V_{n_1}^1=V_{n_1+1}^1=V_{n_1+1}^1=\ldots
		$.
		Since $\mathcal{F}^1_n$ is nonempty for every $n$, we conclude that $V^1$ is nonempty.  As $V^1=V_{n_1}^1$, it must contain $\big(1,F_{n_1}\left(u_1\right)\big)$ for some $F_{n_1}\in\mathcal{F}_{n_1}^1$.  We claim that, for every $n\geq n_1$, there exists $F_n\in\mathcal{F}_n^1$ such that 
		\begin{align}
			F_n\left(u_1\right)=F_{n_1}\left(u_1\right)=:v_1\,.
		\end{align}
		
		To verify the claim above, we observe a trivial fact that $\mathcal{F}_n^1$ is closed under \emph{unit linear combinations}.  That is, if $F_{n,1},F_{n,2},\ldots,F_{n,k}\in \mathcal{F}_n^1$, then 
		$
			\sum_{i=1}^k\,t_i\,F_{n,i}\in \mathcal{F}_n^1
		$
		for every $t_1,t_2,\ldots,t_k\in \amsbb{K}$ with $\sum_{i=1}^k\,t_i=1$.  Therefore, if $\left(1,\tilde{u}_n\right)\in V_n^1$, then 
		$
			\tilde{u}_n=\sum_{i=1}^k\,t_i\,F_{n,i}\left(u_1\right)\,,
		$
		for some $t_1,\ldots,t_k\in\amsbb{K}$ with $\sum_{i=1}^k\,t_i=1$ and for some $F_{n,1},\ldots,F_{n,k}\in\mathcal{F}_n^1$, whence with 
		$
			F_n:=\sum_{i=1}^k\,t_i\,F_{n,i}\in\mathcal{F}_n^1\,,
		$ we have $\tilde{u}_n=F_n\left(u_1\right)$.   In particular, for every positive integer $n$, $v_1$ is in the image of $F_n$ for some $F_n\in\mathcal{F}_n^1$.
		
		Now, suppose that $v_1,v_2,\ldots,v_l$ have been obtained so that, for each $n\in\amsbb{Z}_{>0}$, there exists a map $F_n:N_n\to I_n$ such that 
		\begin{align}
			F_n\circ\varphi_n=f_n\textup{ and }F_n\left(u_i\right)=v_i\textup{ for every }i=1,2,\ldots,l\,.
		\label{eq:injcond}
		\end{align}  
		Let $\mathcal{F}_n^l$ denote the set of all $\Ulie\left(\gfrak'_n\right)$-module homomorphisms $F_n$ that obey (\ref{eq:injcond}).  We proceed as before.  Take $V^{l+1}_n \subseteq \amsbb{K}\times I^{\mu_{l+1}}$ to be the $\amsbb{K}$-span of all vectors of the form $\big(1,F_n\left(u_{l+1}\right)\big)$.  Then, 
		$
			V^{l+1}_1\supseteq V^{l+1}_2\supseteq \ldots$, so that 
		$
			V^{l+1}:=V^{l+1}_{n_{l+1}}=V^{l+1}_{n_{l+1}+1}=\ldots
		$ for some positive integer $n_l$.  Then, $V^{l+1}$ is nonzero and contains $\big(1,F_{n_{l+1}}\left(u_{l+1}\right)\big)$ for some $F_{n_{l+1}}\left(u_{l+1}\right)$.  Then, we set $v_{l+1}$ to be $F_{n_{l+1}}\left(u_{l+1}\right)$.  As before, using the fact that $\mathcal{F}^{l+1}_n$ is closed under unit linear combinations, we conclude that, for every positive integer $n$, there exists $F_n\in \mathcal{F}^{l+1}_n$ for which $F_n\left(u_i\right)=v_i$ for $i=1,2,\ldots,l+1$.

		With known values of $v_1,v_2,\ldots$, we can define $F:N\to I$ via extending the conditions 
		\begin{align}
			F\left(u_i\right)=v_i\textup{ for every }i=1,2,\ldots\,.
		\end{align}  This gives a well defined map as $u_1,u_2,\ldots$ generate $N$.  By the construction, $F\circ \varphi=f$, so that $I$ is injective.
	\end{pf}
\end{thm}
	
	\begin{thm}
		Let $\lambda\in\hfrak^*$ be almost dominant.  Then, there exists an injective hull $\Ilie(\lambda)$ of the simple module $\Llie(\lambda)$.  In particular, if $\lambda$ is dominant, then $\Ilie(\lambda)=\Vlie(\lambda)$.
		\label{thm:almostdominant}
		
		\begin{pf}
			For each positive integer $n$, we write $\Llie_n(\lambda)$ for the simple module in $\bggO^{\gfrak'_n}_{\bfrak_n'}$ with highest weight $\lambda\in\hfrak^*$ as in the proof of Proposition~\ref{prop:blocks}, and denote by $\Ilie_n(\lambda)$ its injective hull $\Ilie\left(\lambda;\gfrak'_n,\bfrak'_n,\hfrak\right)$ in $\bggO^{\gfrak'_n}_{\bfrak_n'}$.
			Similarly, $\Mlie_n(\lambda)$ and $\Vlie_n(\lambda)$ are, respectively, the Verma module $\Mlie\left(\lambda;\gfrak'_n,\bfrak'_n,\hfrak\right)$ and the co-Verma module $\Vlie\left(\lambda;\gfrak_n',\bfrak'_n,\hfrak\right)$ in $\bggO^{\gfrak'_n}_{\bfrak'_n}$ with highest weight $\lambda$.  
			
			We have\footnote{Here, $\ch(M)$ denotes the formal character of an $\hfrak$-weight $\gfrak$-module $M$.}
			$
				\ch\big(\Ilie_n(\lambda)\big)=\sum_{\mu_n\succeq \lambda_n}\,\big\{\Ilie_n(\lambda):\Vlie_n\left(\mu\right)\big\}\,\ch\big(\Vlie_n\left(\mu\right)\big)$.
			Using the finite-dimensional BGG Reciprocity, we have
			\begin{equation}
				\big\{\Ilie_n(\lambda):\Vlie_n\left(\mu\right)\big\}=\big[\Mlie_n\left(\mu\right):\Llie_n\left(\lambda\right)\big]\,.
			\end{equation}
			For each $\mu\succeq \lambda$, there exists $n_\mu\in\amsbb{Z}_{> 0}$ (the existence of $n_\mu$ can be proven via a formal character argument) such that, for all $n\geq n_\mu$, we have
			\begin{equation}
				\big[\Mlie(\mu),\Llie(\lambda)\big]=\big[\Mlie_n\left(\mu\right):\Llie_n\left(\lambda\right)\big]\,.
			\end{equation}
			Because $\lambda$ is almost dominant, there are finitely many $\mu\succeq\lambda$ with $\mu\in W\cdot \lambda$.  Furthermore, the multiplicity $\big\{\Ilie_n(\lambda):\Vlie_n\left(\mu\right)\big\}$  eventually stabilizes at the value $\big[\Mlie(\mu):\Llie(\lambda)\big]<\infty$.  
			
			We have a sequence of embeddings
			$
				\Llie_n\left(\lambda\right)\to\Llie_{n+1}\left(\lambda\right)\to \Ilie_{n+1}(\lambda)$.
			By Proposition~\ref{prop:injrestricted}, $\textup{Res}_{\gfrak_n'}^{\gfrak'_{n+1}}\big(\Ilie_{n+1}(\lambda)\big)$ is injective in $\bggO^{\gfrak'_{n}}_{\bfrak'_{n}}$.  Since $\Ilie_n(\lambda)$ is the injective hull of $\Llie_n\left(\lambda\right)$ in $\bggO^{\gfrak'_n}_{\bfrak'_n}$, there exists an embedding $\Ilie_n(\lambda)\to \Ilie_{n+1}(\lambda)$.
			
			From the work above, we conclude that every weight space of 
			$
				\Ilie(\lambda):=\lim\limits_{\underset{n}{\boldsymbol{\longrightarrow}}}\,\Ilie_n(\lambda)
			$ is finite-dimensional.  This means $\Ilie(\lambda)\in\bar\bggO$, whence $\Ilie(\lambda)$ injective by the previous proposition.  In particular, if $\lambda$ is already dominant, then $\Ilie_n(\lambda)=\Vlie_n\left(\lambda\right)$ for every $n$.  Since the direct limit of $\Vlie_n\left(\lambda\right)$ is just $\Vlie(\lambda)$, the claim follows.
		\end{pf}
	\end{thm}
	
	\begin{thm}
		For a fixed $\lambda\in\hfrak^*$ and $\mu\preceq\lambda$, $\Llie(\mu)$ has an injective hull in $\bar\bggO^{\preceq \lambda}$.
		
		\begin{pf}
			The proof is similar to that of Theorem~\ref{thm:injdirlim}.  We only need to show that the direct limit 
			\begin{align}
				\Ilie^{\preceq\lambda}(\mu):=\lim\limits_{\underset{n}{\boldsymbol{\longrightarrow}}}\,\tcyr_{\preceq\lambda}\Ilie_n\left(\mu\right)
				\end{align}
			 is in $\bbggO^{\preceq\lambda}$, where $\tcyr_{\preceq\lambda}$ also denotes the truncation functor in $\bggO^{\gfrak'_n}_{\bfrak_n'}$ with upper bound $\lambda\in\hfrak^*$.  To this end, we need to verify that $\Ilie^{\preceq\lambda}(\mu)$ has finite-dimensional weight spaces.  
			
			We say that two formal characters $\xi$ and $\zeta$ satisfies $\xi\leq \zeta$ if all coefficients of $e^{\lambda}$ in $\zeta-\xi$ are nonnegative integers.  By studying the formal character of $\tcyr_{\preceq \lambda}\Ilie_n\left(\mu\right)$, it is easy to see that 
			\begin{align}
				\ch\Big(\tcyr_{\preceq\lambda}\Ilie_n\left(\mu\right)\Big) &\leq \sum_{\nu\in\left[\mu,\lambda\right]}\,\big\{\Ilie_n(\mu):\Vlie_n\left(\nu\right)\big\}\,\ch\big(\Vlie_n\left(\nu\right)\big)
				\nonumber\\
				&=\sum_{\nu\in\left[\mu,\lambda\right]}\,\big[\Vlie_n(\nu):\Llie_n\left(\mu\right)\big\}\,\ch\big(\Vlie_n\left(\nu\right)\big)\,.
				\label{eq:truncchar}
			\end{align}
			The right-hand side of (\ref{eq:truncchar}) is bounded as $n\to\infty$.  Therefore, the direct limit $\Ilie^{\preceq\lambda}\left(\mu\right)$ is indeed an object in $\bar\bggO^{\preceq\lambda}$. 
			
			 Since each $\Ilie_n\left(\mu\right)$ is an essential extension of $\Llie_n\left(\mu\right)$ in $\bggO^{\gfrak'_n}_{\bfrak'_n}$, the truncation $\tcyr_{\preceq\lambda}\Ilie_n\left(\mu\right)$ is also an essential extension of $\Llie_n\left(\mu\right)$ in $\left(\bggO^{\gfrak'_n}_{\bfrak'_n}\right)^{\preceq\lambda}$.  Thus, $\Ilie^{\preceq\lambda}\left(\mu\right)$ is indeed the injective hull of $\Llie(\mu)$ in $\bar{\bggO}^{\preceq\lambda}$.
		\end{pf}
	\end{thm}
	
\subsection{BGG Reciprocity}

In this subsection, we shall establish a version of BGG reciprocity for the category $\left(\bbggO^{\gfrak}_{\bfrak}\right)^{\preceq\lambda}$ for a fixed $\lambda\in\hfrak^*$.  To do so, we first note that every object $M$ of $\left(\bbggO^{\gfrak}_{\bfrak}\right)^{\preceq\lambda}$ is countable-dimensional.  Therefore, $M$ can be generated by countably many weight vectors $v_1,v_2,\ldots$.   We set 
$
	M_n:=\sum_{i=1}^n\,\Ulie\left(\gfrak'_n\right)\cdot v_i$.   
The injective hull of $M_n$ in $\left(\bggO^{\gfrak'_n}_{\bfrak'_n}\right)^{\preceq\lambda}$ is denoted by $I_n$.

Since we have a $\gfrak'_n$-module embedding $M_n\to M_{n+1}$ and $I_{n+1}$ is injective as an object of $\left(\bggO^{\gfrak'_n}_{\bfrak'_n}\right)^{\preceq\lambda}$, there exists a $\gfrak'_n$-module embedding $I_n\to I_{n+1}$.  The question is now whether the direct limit $I:=\lim\limits_{\underset{n}{\boldsymbol{\longrightarrow}}}\,I_n$ is in $\left(\bbggO^{\gfrak}_{\bfrak}\right)^{\preceq\lambda}$; that is, we need to check whether the weight spaces of $I$ are finite dimensional.

Let $\mu\in\hfrak^*$ be such that $\mu\preceq \lambda$.  We want to find $\dim_\amsbb{K}\left(I^\mu\right)$.  To do this, we find a bound on $\dim_\amsbb{K}\left(I_n^{\mu}\right)$.  There are at most $\dim_\amsbb{K}\left(M^\mu\right)$ indecomposable direct summands of $I_n$ having $\mu$ as a weight.  We focus on one of such indecomposable direct summands.  It is of the form $\tcyr_{\preceq \lambda}\Ilie_n\left(\xi\right)$ for some $\xi\preceq \mu$ (here, the $\gfrak'_n$-module $\Ilie_n(\xi)$, as well as $\Vlie_n(\nu)$, is as defined in the proof of Theorem~\ref{thm:almostdominant}).  

The contribution to the weight space with weight $\mu$ of $\tcyr_{\preceq \lambda}\Ilie_n\left(\xi\right)$ can only come from its co-Verma subquotients $\Vlie_n\left(\nu\right)$ with $\mu\preceq \nu\preceq\lambda$.  Thus, we have an upper bound
\begin{align}
	\dim_\amsbb{K}\left(I_n^{\mu}\right)\leq \sum_{\xi}\,\sum_{\nu\in\left[\mu,\lambda\right]}\,\Big\{\tcyr_{\preceq \lambda}\Ilie_n\left(\xi\right):\Vlie_n\left(\nu\right)\Big\}\,\dim_\amsbb{K}\left(\Vlie_n\left(\nu\right)^{\mu}\right)\,,
\end{align}
where $\xi$ runs over possible weights such that $\tcyr_{\preceq \lambda}\Ilie\left(\xi\right)$ is an indecomposable direct summand of $I_n$ with $\mu$ as a weight, and $[\mu,\lambda]$ denotes the set $\big\{\nu\in\hfrak^*\suchthat{}\mu\preceq\nu\preceq\lambda\big\}$.
By the BGG reciprocity, we have 
\begin{align}
	\dim_\amsbb{K}\left(I_n^{\mu}\right)&\leq \sum_{\xi}\,\sum_{\nu\in\left[\mu,\lambda\right]}\,\Big[\Mlie_n\left(\nu\right):\Llie_n\left(\xi\right)\Big]\,\dim_\amsbb{K}\left(\Mlie_n\left(\nu\right)^{\mu}\right)
	\nonumber \\
	&\leq \dim_\amsbb{K}\left(M^\mu\right)\,\sum_{\nu\in\left[\mu,\lambda\right]}\,A_n(\nu)\,\dim_\amsbb{K}\left(\Mlie\left(\nu\right)^{\mu}\right)\,,
	\label{eq:injbound}
\end{align}
where $A_n(\nu)$ is the maximum possible value of $\Big[\Mlie_n\left(\nu\right):\Llie_n\left(\xi\right)\Big]$ with $\xi\preceq \mu$.

We are now ready to prove the proposition below.

\begin{prop}
	If $\lambda$ is an almost antidominant weight, then the truncated category $\left(\bbggO^{\gfrak}_{\bfrak}\right)^{\preceq\lambda}$ has enough injectives (and so, $\left(\bbggO^{\gfrak}_{\bfrak}\right)^{\preceq\lambda}$ has enough projectives as well).
	
	\begin{pf}
		If $\lambda$ is almost antidominant, then $\mu$ is also almost antidominant.  Therefore, there are finitely many weights $\xi\in\hfrak^*$ such that $\xi\preceq \mu$.  Thus, if $A$ denotes the maximum of $\Big[\Mlie(\nu):\Llie(\xi)\Big]$ with $\xi\preceq \mu$ and $\nu\in[\mu,\lambda]$, we have from (\ref{eq:injbound}) that
\begin{align}
	\dim_\amsbb{K}\left(I_n^{\mu}\right) &\leq \dim_\amsbb{K}\left(M^\mu\right)\,\sum_{\nu\in\left[\mu,\lambda\right]}\,A_n(\nu)\,\dim_\amsbb{K}\left(\Mlie\left(\nu\right)^{\mu}\right)
	\nonumber\\
	&\leq A\,\dim_\amsbb{K}\left(M^\mu\right)\,\sum_{\nu\in\left[\mu,\lambda\right]}\,\dim_\amsbb{K}\left(\Mlie\left(\nu\right)^{\mu}\right)<\infty\,,
\end{align}
whenever $n$ is large enough.  Ergo, there exists a universal bound for the dimension of the weight space $I_n^{\mu}$ for all (sufficiently large) $n$.  That is, $\dim_\amsbb{K}(I^\mu)<\infty$ and the claims follows immediately.
	\end{pf}
\end{prop}

Now, we want to show that, for any $\lambda\in\hfrak^*$ and $\mu\preceq\lambda$, the injective hull $\Ilie^{\preceq\lambda}(\mu)$ has a co-standard filtration.  We recall from the finite-dimensional theory that there exists a co-Verma filtration
\begin{align}
	0=F_n^0\subsetneq F_n^1\subsetneq \ldots \subsetneq F^{t_n-1}_n\subsetneq F^{t_n}_n=\Ilie_n^{\preceq \lambda}\left(\mu\right)\,,
\end{align}
where $\Ilie_n^{\preceq\lambda}\left(\mu\right)$ denotes the $\gfrak'_n$-module $\tcyr_{\preceq\lambda}\,\Ilie_n(\mu)$.  Since the highest weights of $F_n^{i}/F^{i-1}_n$ are in the interval $\left[\mu,\lambda\right]$, we have by BGG reciprocity that
\begin{align}	
	t_n \leq \sum_{\nu\in\left[\mu,\lambda\right]}\,\Big[\Mlie_n\left(\nu\right):\Llie_n\left(\mu\right)\Big] =  \sum_{\nu\in\left[\mu,\lambda\right]}\,\Big[\Mlie\left(\nu\right):\Llie\left(\mu\right)\Big]<\infty
\end{align}
for sufficiently large $n$.  Thus, there exists a sequence $\left(n_k\right)_{k=1}^\infty$ of positive integers such that 
\begin{align}
	n_1<n_2<n_3\ldots
\textup{ and }
	t:=t_{n_1}=t_{n_2} =t_{n_3}=\ldots\,.
\end{align}   
Furthermore, as there are only finitely many weights $\upsilon \in W[\lambda]\cdot \lambda$ with $\mu\preceq\upsilon\preceq \lambda$, we may assume without loss of generality (due to the Pigeonhole Principle) that the highest weight $\xi_{n_k}[i]$ of $F^i_{n_{k+1}}/F^{i-1}_{n_{k+1}}$ is the same as the highest weight of $F^{i}_{n_k}/F^{i-1}_{n_k}$ for every $i$ and $k$.  Denote by $\xi[i]\in\hfrak^*$ the common weight $\xi_{n_1}[i],\xi_{n_2}[i],\xi_{n_3}[i],\ldots$.  We need the following lemma.

\begin{lem}
	Let $m,n\in\amsbb{Z}_{>0}$ be such that $m\leq n$.  The Verma module $\Mlie_m\left(\lambda\right)$ is a direct summand of the Verma module $\Mlie_n\left(\lambda\right)$, viewed as a $\gfrak'_m$-module.  Consequently, the co-Verma module $\Vlie_m\left(\lambda\right)$ is also a direct summand of the co-Verma module $\Vlie_n\left(\lambda\right)$, viewed as a $\gfrak_m'$-module.
	\label{lem:covermaemb}
	
	\begin{pf}
		Let $v$ be a highest-weight vector of $\Mlie_n\left(\lambda\right)$.  Let 
		$
			\left\{x_{\pm\alpha}\,\boldsymbol{|}\,\alpha\in\Delta^+\right\}\cup\left\{h_\beta\,\boldsymbol{|}\,\beta\in\Sigma^+\right\}$
		be a Chevalley basis of $\gfrak$.  Then, 
		\begin{align}
			\textup{Res}^{\gfrak_n'}_{\gfrak_m'}\Mlie_n\left(\lambda\right)=\Ulie\left(\gfrak_n'\right)\cdot v= \Big(\Ulie\left(\gfrak'_m\right)\cdot v \Big)\oplus \left(\sum_{\alpha}\,\Ulie\left(\gfrak'_n\right)\cdot \left(x_{-\alpha}\cdot v\right)\right)\,,
		\end{align}
		where $\alpha$ runs over $\bfrak_n'$-positive roots of $\gfrak'_n$ which are not roots of $\gfrak'_m$, is a direct sum decomposition of $\Mlie_n\left(\lambda\right)$ as a $\gfrak'_m$-module with a direct summand 
		$
			\Ulie\left(\gfrak'_m\right)\cdot v\cong \Mlie_m\left(\lambda\right)
		$.  To prove the co-Verma version, we only need to apply the duality functor.
	\end{pf}
\end{lem}

Clearly, $F^1_{n_k}$ is a co-Verma submodule of $I_k:=\Ilie^{\preceq \lambda}_{n_k}\left(\mu\right)$ containing the socle of $I_k$.  By the lemma below, each $F^1_{n_k}$ is unique as it contains the simple module 
$
\Ulie\left(\gfrak_{n_k}\right)\cdot v[1]\cong\Llie\big(\xi[1]\big)
$, where $v[1]$ is a singular vector of the socle of $\Ilie^{\preceq\lambda}(\mu)$.   Then, using Lemma~\ref{lem:covermaemb}, we have a embeddings $F^1_{n_k}\to F^1_{n_{k+1}}$, whose direct limit is simply the co-Verma module $
	F^1\cong\Vlie\big(\xi[1]\big)$,
where $\xi[1]$ is clearly equal to $\mu$. 

\begin{lem}
	Let $n\in\amsbb{Z}_{>0}$ and $M_n\in\bggO^{\gfrak'_n}_{\bfrak'_n}$.  Suppose that a simple module $L_n\in\bggO^{\gfrak'_n}_{\bfrak'_n}$ is a submodule of $M_n$.  Then, $M_n$ has at most one co-Verma submodule $V_n$ such that $L_n\subseteq V_n\subseteq M_n$.
	\label{lem:uniquecoverma}
	
	\begin{pf}
		Suppose that $M$ has two co-Verma submodules $V_n$ and $V'_n$ with $L_n\subseteq V_n$ and $L_n\subseteq V'_n$.  Take 
		$N_n:=V_n+V_n'$. Then, $N_n$ is indecomposable (as $V_n$ and $V'_n$ are both indecomposable with $V_n\cap V'_n\supseteq L_n\supsetneq 0$).  Hence, we have a short exact sequence
		$
			0 \to V_n \to N_n \to N_n/V_n \to 0$.
		Dualizing this exact sequence yields
		\begin{align}
			0\to \left(N_n/V_n\right)^\vee \to N_n^\vee \to V_n^\vee \to 0\,.
		\end{align}
		By Proposition~\ref{prop:extintro}, we see that this exact sequence must split.  As $N_n^\vee$ is indecomposable, we conclude that $V_n^\vee=0$ or $\left(N_n/V_n\right)^\vee=0$.  Since $V_n\neq 0$, we must have $N_n/V_n=0$, which leads to $V'_n=V_n$.
	\end{pf}
\end{lem}

Suppose now that, for some positive integer $l<t$, the submodules $0=F^0$, $F^1$, $F^2$, $\ldots$, $F^l$ of $\Ilie^{\preceq \lambda}(\mu)$ have been determined with the property that $F^i$ is the direct limit $\lim\limits_{\longrightarrow}\,F_{n_k}^i$, where the $\gfrak'_{n_k}$-modules $F_{n_k}^i$ are submodules of $I_k$ satisfying the following properties:
\begin{enumerate}[(i)]
\item $0=F^0_{n_k}\subsetneq F^1_{n_k} \subsetneq \ldots \subsetneq F^l_{n_k}\subsetneq I_k$,
	\item $F^{i}_{n_k}/F^{i-1}_{n_k}\cong\Vlie\big(\xi[i]\big)$ for every $i=1,2,\ldots,l$.
\end{enumerate}Then, we proceed by looking at the quotient $I_k/F^{l}_{n_k}$.  Identify each $u+F^l_{n_k}\in I_k/F^l_{n_k}$ as an element of $I/F^l$ via
$
	u+F^l_{n_k}\mapsto u+F^l\in I/F^l$
(making $I_k/F_{n_k}^l$ a $\gfrak'_{n_k}$-submodule of $I/F^l$).  We have an embedding $
	\Vlie\big(\xi[l+1]\big) \to I_k/F^l_{n_k}$
for each $k$.  Let $V^{l+1}_k$ be the $\amsbb{K}$-span of all vectors $v\in I_k/F^{l}_{n_k} \subseteq I/F^{l}$ such that $v$ is the image of a singular vector under an embedding $\Vlie\big(\xi[l+1]\big)\to I_k/F^{l+1}_{n_k}$.  Hence, $V^{l+1}_k$ is a nonzero subspace of $\left(I/F^l\right)^{\xi[l+1]}$ and $V^{l+1}_k\supseteq V^{l+1}_{k+1}$ for every $k$.  Because 
\begin{align}
	\dim_\amsbb{K}\left(V_k^{l+1}\right)\leq \dim_\amsbb{K}\left(\left(I/F^l\right)^{\xi[l+1]}\right)<\infty\,,
\end{align} 
there exists $
	v[l+1]+F^l\in \bigcap_{k\in\amsbb{Z}_{>0}}\,V^{l+1}_k$ which is a nonzero element of $\left(I/F^l\right)^{\xi[l+1]}$. 

Now, by Lemma~\ref{lem:uniquecoverma}, we can show that there is a unique co-Verma submodule $\bar{F}^{l+1}_{n_k}$ of $I_k/F_{n_k}^l$ containing the simple submodule $\Ulie\left(\gfrak'_{n_k}\right)\cdot\left(v[l+1]+F^l\right)$.  Then, the direct limit $\bar{F}^{l+1}$ of $\bar{F}^{l+1}_{n_k}$ must be a co-Verma module of highest weight $\xi[l+1]$.  Let $F^{l+1}$ be the preimage of $\bar{F}^{l+1}$ under the quotient map $I\to I/F^l$.   Then, by induction, we have found a filtration
\begin{align}
	0=F^0\subsetneq F^1\subsetneq F^2\subsetneq \ldots\subsetneq F^{t-1}\subsetneq F^t=\Ilie^{\preceq\lambda}(\mu)
	\label{eq:covermafiltind}
\end{align}
of $\Ilie^{\preceq \lambda}(\mu)$ such that each successive quotient $F^l/F^{l-1}$ is isomorphic to the co-Verma module $\Vlie\big(\xi[l]\big)$.   It can be easily seen that the number of times a co-Verma module $\Vlie(\nu)$ appears as a successive quotient $F^l/F^{l-1}$ in (\ref{eq:covermafiltind}) is independent on the choice of the co-Verma filtration.   We use the notation $\Big\{\Ilie^{\preceq\lambda}(\mu):\Vlie(\nu)\Big\}$ for the number of times that $\Vlie(\nu)$ appears as a successive quotient in (\ref{eq:covermafiltind}). 

Let $\Plie^{\preceq \lambda}(\mu)$ denote $\tcyr^\vee_{\preceq \lambda}\Plie(\mu)=\Big(\tcyr_{\preceq\lambda}\Ilie(\mu)\Big)^\vee$.  Then, by applying duality on the co-Verma filtration (\ref{eq:covermafiltind}), $\Plie^{\preceq \lambda}(\mu)$ has a Verma filtration
\begin{align}
	0=T^0\subsetneq T^1\subsetneq T^2\subsetneq \ldots\subsetneq T^{t-1}\subsetneq T^t=\Plie^{\preceq\lambda}(\mu)\,,
	\label{eq:vermafiltind}
\end{align}
where each successive quotient $T^l/T^{l-1}$ is isomorphic to the Verma module $\Mlie\big(\xi[t+1-l]\big)$.  In particular, 
$T^t/T^{t-1}\cong \Mlie(\mu)$.  The number of times that $\Mlie(\nu)$ appears as a successive quotient $T^l/T^{l-1}$ in (\ref{eq:vermafiltind}) is also well defined, and is denoted by $\Big\{\Plie^{\preceq\lambda}(\mu):\Mlie(\nu)\Big\}$.

\begin{prop}
	For every $\lambda,\mu\in\hfrak^*$ with $\mu\preceq \lambda$, the injective object $\Ilie^{\preceq\lambda}(\mu)$ has a finite filtration with successive quotients isomorphic to co-Verma modules.  Furthermore, we have BGG reciprocity:
	\begin{align}
		\Big\{\Plie^{\preceq\lambda}(\mu):\Mlie(\nu)\Big\}=\Big\{\Ilie^{\preceq\lambda}(\mu):\Vlie(\nu)\Big\}=\Big[\Mlie(\nu):\Llie(\mu)\Big]=\Big[\Vlie(\nu):\Llie(\mu)\Big]\,,
	\end{align}
	for all $\nu\in[\mu,\lambda]$.  
	\label{prop:covermafiltration}
\end{prop}

Finally, we note that, if $\lambda$ is not almost antidominant, then $\Mlie(\lambda)$ is of infinite length and cannot be written as a union of subobjects of finite length.  This is because every submodule $M$ of $\Mlie(\lambda)$ has a singular vector $v\neq 0$.  The submodule $N$ of $M$ generated by $v$ is then a Verma module with highest weight $\mu\preceq\lambda$, which is not almost antidominant.  Ergo, $N$ is of infinite length, and so is $M$.  Thus, $\Mlie(\lambda)$ has no submodules of finite length.  In particular, this implies that $\Mlie(\lambda)$ has trivial socle.

The argument above shows that $\left(\bbggO^{\gfrak}_{\bfrak}\right)^{\preceq \lambda}$ is not locally artinian, whence this category is not a highest-weight category in the sense of \cite{hwcat}.  That is, $\left(\bbggO^{\gfrak}_{\bfrak}\right)^{\preceq \lambda}$ is not a highest-weight category.  Combining this observation with the fact that $\left(\bbggO^{\gfrak}_{\bfrak}\right)^{\preceq \lambda}$ has enough injectives when $\lambda$ is almost antidominant, we conclude the following theorem.

\begin{thm}
	The category $\left(\bbggO^\gfrak_\bfrak\right)^{\preceq \lambda}$ is a highest-weight category if and only if $\lambda$ is an almost antidominant weight.  
\end{thm}

\begin{openq}
For a weight $\lambda\in\hfrak^*$ which is not almost antidominant, does the category $\left(\bbggO^{\gfrak}_{\bfrak}\right)^{\preceq\lambda}$ have enough injectives?
\end{openq}

\begin{define}
	Let $\mathcal{C}$ be an abelian category with an abelian subcategory $\tilde{\mathcal{C}}$.  An object $I\in\mathcal{C}$ is \emph{injective relative to $\tilde{\mathcal{C}}$} if, for any two objects $X,Y\in\tilde{\mathcal{C}}$ and any monomorphism $f\in \Hom_{\tilde{\mathcal{C}}}(X,Y)$, every morphism $g\in \Hom_{\mathcal{C}}(X,I)$ factors through $f$, i.e., there exists $\varphi\in\Hom_{\mathcal{C}}(Y,I)$ such that 
	$
		g=\varphi\circ f$.
\end{define}

\begin{thm}
	Let $R$ be a ring.  Suppose that $\mathcal{C}$ and $\tilde{\mathcal{C}}$ are abelian subcategories of the category of left $R$-modules with $\tilde{\mathcal{C}}$ being a subcategory of $\mathcal{C}$.  If $M\in\tilde{\mathcal{C}}$ has an injective hull $I$ in $\tilde{\mathcal{C}}$, then for each object $J\in\mathcal{C}$ which is injective relative to $\tilde{\mathcal{C}}$, any embedding $\iota\in\Hom_{\mathcal{C}}(M,J)$ induces an embedding $\varphi\in\Hom_\mathcal{C}(I,J)$.
	\label{thm:relinj}
	
	\begin{pf}
		We have an exact sequence $0\to M\to I$ of objects and morphisms in $\tilde{\mathcal{C}}$ and a homomorphism $\iota\in\Hom_\mathcal{C}(M,J)$.  As $J$ is injective relative to $\tilde{\mathcal{C}}$, there exists a map $\varphi\in\Hom_{\mathcal{C}}(I,J)$ such that the diagram below commutes:
		\begin{equation}
		\begin{tikzcd}
		0\arrow{r}& M \arrow{r}{\subseteq} \arrow{d}{\iota} & I\arrow{dl}{\varphi}\\
		&J \,.&
		\end{tikzcd}
		\label{eq:relinj}
		\end{equation}
		We claim that $\varphi:I\to J$ is an embedding. 
		
		Let $K:=\ker(\varphi)$.  Since $\varphi|_M=\iota$ due to commutativity of (\ref{eq:relinj}) and $\iota$ is an embedding, we must have 
		$
			K\cap M=\ker\left(\varphi|_M\right)=\ker(\iota)=0
		$.  Because $I$ is an essential extension of $M$, the condition $K\cap M=0$ implies that $K=0$.  Therefore, $\varphi$ is injective.
	\end{pf}
\end{thm}

	\begin{thm}
		For $\lambda\in\hfrak^*$, the simple module $\Llie(\lambda)$ has an injective hull and a projective cover in $\bbggO$ if and only if $\lambda$ is almost dominant.   In particular, this implies that $\bbggO$ does  not have enough injectives, and therefore, $\bbggO$ is not a highest-weight category.
		
		\begin{pf}
			If $\lambda$ is almost dominant, then Theorem~\ref{thm:almostdominant} shows that $\Llie(\lambda)$ has an injective hull in $\bbggO$, and by duality, it has also a projective cover.  To prove the converse, we suppose on the contrary that $\lambda$ is not almost dominant but $\Llie(\lambda)$ has an injective hull $I$ in $\bbggO$.  
			
			As $\lambda$ is not almost dominant, there exists a sequence of weights $\left(\lambda_i\right)_{i=0}^\infty$ with $\lambda_i\in W[\lambda]\cdot \lambda$ and 
			\begin{align}
				\lambda =\lambda_0 \prec \lambda_1 \prec \lambda_2 \prec \ldots\,.
			\end{align}
			For simplicity, let $I_i$ denote $\Ilie^{\preceq \lambda_i}\left(\lambda\right)$ for $i=0,1,2,\ldots$.  It is clear that $I$ is injective relative to $\bbggO^{\preceq\lambda_i}$ for each $i$.  By Theorem~\ref{thm:relinj}, there exists an embedding of $I_i$ into $I$.
			
			Now, using Proposition~\ref{prop:covermafiltration}, we know that each $I_i$ has a co-Verma filtration 
			\begin{align}
				\Vlie(\lambda)=F_i[0]\subsetneq F_i[1] \subsetneq \ldots \subsetneq F_i\left[k_i\right]=I_i\,.
				\label{eq:covermafilt}
			\end{align}
			Furthermore, as $I_i \subsetneq I_{i+1}$, we have 
			$k_0< k_1< k_2< \ldots$.  For every $j=1,2,\ldots,k_i$, the successive quotient $F_i[j]/F_i[j-1]$ is isomorphic to the co-Verma module $\Vlie\left(\mu_i[j]\right)$ for some $\mu_i[j]\in W[\lambda]\cdot\lambda$ with $\mu_j\succeq \lambda$.  This implies 
			$\dim_\amsbb{K}\left(\big(F_i[j]/F_i[j-1]\big)^\lambda \right)\geq 1$.  Ergo,
			\begin{align}
				\dim_\amsbb{K}\left(I^\lambda\right) \geq \dim_\amsbb{K}\left(I_i^\lambda\right) \geq \sum_{j=1}^{k_i}\,\dim_\amsbb{K}\left(\big(F_i[j]/F_i[j-1]\big)^\lambda \right)\geq k_i
			\end{align}
			for every $i=0,1,2,\ldots$.  As $\lim\limits_{i\tendsto\infty}\,k_i=\infty$, we conclude that 
			$\dim_\amsbb{K}\left(I^\lambda\right)=\infty$, which is absurd.  Hence, $\Llie(\lambda)$ does not have an injective hull in $\bbggO$.  
			
			Using duality, we also conclude that $\Llie(\lambda)$ does not have a projective cover in $\bbggO$.  The theorem follows.
		\end{pf}
	\end{thm}
%

The theorem above is the reason why we need to truncate the category $\bbggO$.  With truncation, every simple object has an injective hull and a projective cover, and because of that, a version of BGG reciprocity holds.


\begin{thebibliography}{9}
\bibliography{}

\addcontentsline{toc}{section}{References}

\bibitem{BS}
	A. Baranov and H. Strade, \href{http://www.sciencedirect.com/science/article/pii/S0021869302000790}{\emph{Finitary Lie Algebras}}, Journal of Algebra, {\bf254} (2002), 173--211.
	


%

\bibitem{BK}
	J. L. Brylinski and M. Kashiwara, \href{https://eudml.org/doc/142827}{\emph{Kazhdan-Lusztig Conjecture and Holonomic Systems}}, Inventiones Mathematicae, {\bf64} (1981), 387--410.
	
\bibitem{Penkov-Chirvasitu}
	A. Chirvasitu and I. Penkov, \href{http://math.jacobs-university.de/penkov/papers/m16.pdf}{\emph{Categories of Tensor Representations for Mackey Lie Algebras}}, preprint 2017.


\bibitem{hwcat}
	E. Cline, B. Parshall, and L. Scott, \href{http://u.math.biu.ac.il/~margolis/Representation%20Theory%20Seminar/Highest%20Weight%20Categories.pdf}{\emph{Finite-Dimensional Algebras and Highest-Weight Categories}}, Journal f\"ur die reine und angewandte Mathematik, {\bf391} (1988), 85--99.

\bibitem{dancohen}
	E. Dan-Cohen, \href{http://books.google.de/books?id=RjjebpvxliMC}{\emph{Structure of Root-Reductive Lie Algebras}}, Ph.D. Dissertation, University of California, Berkeley, (2002).
	
\bibitem{BorelDC}
	E. Dan-Cohen, \href{https://arxiv.org/abs/0708.1267}{\emph{Borel Subalgebras of Root-Reductive Lie Algebras}}, Journal of Lie Theory, {\bf18} (2002), 215--241.

\bibitem{DPS}
	E. Dan-Cohen, I. Penkov, and N. Snyder, \href{http://math.jacobs-university.de/penkov/papers/dps.pdf}{\emph{Cartan Subalgebras of Root-Reductive Lie Algebras}}, Journal of Algebra, {\bf308} (2007), 583--661.	

	
	
\bibitem{DP1999}
	I. Dimitrov and I. Penkov, \href{http://www.mast.queensu.ca/~dimitrov/DP3.pdf}{\emph{Weight Modules of Direct Limit Lie Algebras}}, International Mathematics Research Notices, {\bf5} (1999), 223--249.
	
	
	

\bibitem{KSRA}
	A. Facchini, \href{http://www.springer.com/cda/content/document/cda_downloaddocument/9783034803021-c1.pdf?SGWID=0-0-45-1316849-p174271584}{\emph{Module Theory: Endomorphism Rings and Direct Sum Decompositions in Some Classes of Modules}}, Birkh\"{a}user, New York, (2010).
	
\bibitem{FF}
	B. Feigin and E. Frenkel, \href{http://www.worldscientific.com/doi/abs/10.1142/S0217751X92003781}{\emph{Affine Kac-Moody Algebras at the Critical Level and Gelfand-Dikii Algebras}}, International Journal of Modern Physics A, {\bf07} (1992), 197--215. 


	
\bibitem{humphreys}
	J. E. Humphreys, \href{http://books.google.de/books/about/Introduction_to_Lie_algebras_and_represe.html?id=mSD4DvUFa6QC}{\emph{Introduction to Lie Algebras and Representation Theory}}, Springer, New York, (1970).
	

\bibitem{Coxeter}
	J. E. Humphreys, \href{https://books.google.de/books/about/Reflection_Groups_and_Coxeter_Groups.html?id=ODfjmOeNLMUC&redir_esc=y}{\emph{Reflection Groups and Coxeter Groups}}, Cambridge University Press, New York, (1992).

\bibitem{bggo}
	J. E. Humphreys, \href{http://books.google.de/books/about/Representations_of_Semisimple_Lie_Algebr.html?id=8GCP4Ng6risC&redir_esc=y}{\emph{Representations of Semisimple Lie Algebras in the BGG Category $\pazocal{O}$}}, American Mathematical Society, Providence, (2000).
	
\bibitem{kac}
	V. G. Kac, \href{http://books.google.de/books/about/Infinite_Dimensional_Lie_Algebras.html?id=kuEjSb9teJwC&redir_esc=y}{\emph{Infinite-Dimensional Lie Algebras}}, 3rd Edition, Cambridge University Press, New York, (1994).
	
\bibitem{KL}
	D. Kazhdan and G. Lusztig, \href{https://link.springer.com/article/10.1007%2FBF01390031?LI=true}{\emph{Representations of Coxeter Groups and Hecke Algebras}}, Inventiones Mathematicae, {\bf53} (1979), 165--184.
	
\bibitem{NP}
	K.-H. Neeb and I. Penkov, \href{https://www.algeo.math.fau.de/fileadmin/algeo/users/neeb/pdf-artikel/carsubalg.pdf}{\emph{Cartan Subalgebras of $\gllie_\infty$}}, Canadian Mathematical Bulletin, {\bf46} (2003), 597--616.


\bibitem{penkov-dimitrov}
	I. Penkov and I. Dimitrov, \href{http://math.jacobs-university.de/penkov/papers/borel_subalg_dimitrov_penkov.pdf}{\emph{Borel Subalgebras of $\gllie(\infty)$}}, Resenhas IME-USP, {\bf 6} (2007), 111--119.
	
\bibitem{penkov-strade}
	I. Penkov and H. Strade, \href{http://link.springer.com/article/10.1007\%2Fs00013-003-0791-3}{\emph{Locally Finite Lie Algebras with Root Decomposition}}, Archiv der Mathematik, {\bf 80} (2003), 478--485.
	

\bibitem{RCW}
	A. Rocha-Caridi and N. R. Wallach, \href{http://link.springer.com/article/10.1007%2FBF01318901?LI=true}{\emph{Projective modules over graded Lie algebras. I}}, Mathematische Zeitschrift, {\bf180} (1982), 151--177.


\bibitem{SVV}
	P. Shan, M. Varagnolo, and E. Vasserot, \href{https://arxiv.org/pdf/1107.0146.pdf}{\emph{Koszul Duality of Affine Kac-Moody Algebras and Cyclotomic Rational Double Affine Hecke Algebras}}, Advances in Mathematics, {\bf262} (2014), 370--435.
	

	

\end{thebibliography}


\Addresses

\end{document}